%
%
%
%
%
%
%
%
%

%

\def\State#1 {\statementtag#1 }

\catcode`\^^J=10
\magnification=\mag
\documentstyle{amsppt}
\pagewidth{6.00 truein}
\pageheight{8.00 truein}
\hcorrection{.4 truein}
\vcorrection{.25 truein}

\topskip= 28pt

\TagsOnRight

%

\def\ikedocument{
\ifbibmakemode\immediate\openout1= \jobname.bib\fi
\ifmultisection\immediate\openout2=\jobname.eqn\fi}

%

\output={\plainoutput}
\headline={\hfil}
\footline={\hss\bf\folio\hss}

%
%
%
%
%

\newif\ifproofmode    \proofmodefalse
\newif\ifbibmakemode  \bibmakemodefalse
\newif\ifbibcitemode  \bibcitemodefalse
\newif\ifmultisection \multisectionfalse
\newif\ifmulteq        \multeqtrue

%

\newcount\lastsectpageno
\lastsectpageno=0
\define\sectpagenocs#1{sect\expandafter\romannumeral#1pageno}

\define\initializesection#1{\sectionnumber=#1%
\equationnumber=1 \statementnumber=0\exercisenumber=0
\ifmultisection
\global\advance\lastsectpageno by 1
\immediate\write2{\noexpand\expandafter\def%
\noexpand\csname\sectpagenocs{\the\sectionnumber}\noexpand\endcsname
{\the\pageno}}
\fi}

\define\forwardsectpageno{\immediate\write2{\lastsectpageno=\the\pageno}}

\define\tocpageno#1{
\expandafter\csname\sectpagenocs{#1}\expandafter\endcsname}

%

\newcount\equationnumber
\newcount\eqnolet
\newcount\sectionnumber
\newcount\statementnumber
\newcount\exercisenumber
\def\strutdepth{\dp\strutbox}

\def\margintagleft#1{\strut\vadjust{\kern-\strutdepth
{\vtop to \strutdepth{\baselineskip\strutdepth\vss\llap{\sevenrm
#1\quad}\null}}}}

\def\ifundefined#1{\expandafter\ifx\csname#1\endcsname\relax}%

\def\eqcs#1{s\expandafter\romannumeral\the\sectionnumber%
eq\romannumeral#1}

\def\writeeqcs#1{\ifmultisection%
\immediate\write2{\noexpand\expandafter\def%
\noexpand\csname\eqcs#1\noexpand\endcsname{\the\equationnumber}}\fi%
}

\def\makeeqcs#1{\expandafter\xdef\csname\eqcs#1\endcsname%
{\the\equationnumber}%
}

\def\equationtag#1#2{
\ifundefined{\eqcs#1}
\else
\message{Tag \the\sectionnumber.#1.#2 already exists}
\fi%
\tag"(\the\sectionnumber.\the\equationnumber#2)%
\ifproofmode\rlap{\quad\sevenrm#1}\fi%
\writeeqcs{#1}\makeeqcs{#1}%
\global\eqnolet=\equationnumber%
\global\advance\equationnumber by 1%
\global\def\eqrefnumber{#1}"%
}

\def\Tag#1 {\equationtag{#1}{\empty}}

\define\Tagletter#1 {\tag\the\sectionnumber.\the\eqnolet#1}

\define\equationlabel#1#2#3{%
{\xdef\cmmm{\csname
 s\romannumeral#1eq\romannumeral#2\endcsname}\expandafter\ifx
 \csname s\romannumeral#1eq\romannumeral#2\endcsname
\relax\message{Equation #1.#2.#3 not defined ...}\fi(#1.\cmmm#3)}}

\define\eq#1{\equationlabel{\the\sectionnumber}{#1}{}}

%

\define\statementtag#1 {%
\expandafter\ifx\csname %
s\expandafter\romannumeral\the\sectionnumber %
stat\romannumeral#1\endcsname\relax\else%
\message{Statementtag \the\sectionnumber.#1 already exists ...}\fi%
\ifproofmode%
\margintagleft{#1}\fi\global\advance\statementnumber by 1%
\expandafter\xdef\csname %
s\expandafter\romannumeral\the\sectionnumber %
 stat\romannumeral#1\endcsname{\the\statementnumber}\ifmultisection %
\immediate\write2{\noexpand\expandafter\def\noexpand\csname %
s\expandafter\romannumeral\the\sectionnumber%
stat\romannumeral#1%
\noexpand\endcsname{\the\statementnumber}}%
\fi\the\sectionnumber.\the\statementnumber}

\define\statementlabel#1#2{\xdef\cmmm{\csname
s\romannumeral#1stat\romannumeral#2\endcsname}\expandafter\ifx
 \csname s\romannumeral#1stat\romannumeral#2\endcsname\relax\message{Statement
 #1.#2 not defined ...}\fi#1.\cmmm}

\define\statement#1{\statementlabel{\the\sectionnumber}{#1}}
\define\st#1{\statementlabel{\the\sectionnumber}{#1}}


\define\exercisetag#1 {%
\expandafter\ifx\csname %
s\expandafter\romannumeral\the\sectionnumber %
exer\romannumeral#1\endcsname\relax\else%
\message{exercisetag \the\sectionnumber.#1 already exists ...}\fi%
\ifproofmode%
\margintagleft{#1}\fi\global\advance\exercisenumber by 1%
\expandafter\xdef\csname %
s\expandafter\romannumeral\the\sectionnumber %
 exer\romannumeral#1\endcsname{\the\exercisenumber}\ifmultisection %
\immediate\write2{\noexpand\expandafter\def\noexpand\csname %
s\expandafter\romannumeral\the\sectionnumber%
exer\romannumeral#1%
\noexpand\endcsname{\the\exercisenumber}}%
\fi{\bf\the\sectionnumber.\the\exercisenumber.\quad }}

\define\exerciselabel#1#2{\xdef\cmmm{\csname
s\romannumeral#1exer\romannumeral#2\endcsname}\expandafter\ifx
 \csname s\romannumeral#1stat\romannumeral#2\endcsname\relax\message{exercise
 #1.#2 not defined ...}\fi#1.\cmmm}
\define\Exer#1 {\noindent\exercisetag{#1} }

\def\ex#1{\exerciselabel{\the\sectionnumber}{#1}}
\def\Sol#1{\noindent {\bf \exerciselabel{\the\sectionnumber}{#1}\quad}}

%
%
%
%
%

\def\bibyear#1:#2 {\edef\reftag{#1\romannumeral#2}}

\font\slr=cmsl10

\def\cite#1{\catcode`\- =11
\ifbibmakemode\immediate\write1{#1}[{\bf00}]%
\fi%
\ifbibcitemode%
\bibyear#1  %
\expandafter\ifx\csname\reftag bibno\endcsname\relax{\message{#1 not in
bibfile}[{\bf00}]}%
\else%
[{\slr \csname\reftag bibno\endcsname}\hbox{\kern 1pt}]%
\fi\fi
}

\def\referencetag#1:#2 {\edef\reftag{#1\romannumeral#2}%
\ifbibmakemode%
\immediate\write3{\noexpand\def\csname \reftag bibno\endcsname
{\the\refnumb}}
\fi}

\newcount\refnumb
\refnumb=0
\def\InitializeRef{
\ifbibmakemode\immediate\openout3=\jobname.ref\fi
\NoBlackBoxes
\medskip}

%

\define\today{\ifcase\month\or January \or February \or March
\or April  \or May \or June \or July \or August \or September
\or October \or November \or December \fi\space
\oldnos{\number\day}, \oldnos{\number\year}}

\def\hook{\mathbin{\raise2.5pt\hbox{\hbox{{\vbox{\hrule height.4pt
width6pt depth0pt}}}\vrule height3pt width.4pt depth0pt}\,}}

%

\define\Ex#1 {\par\medpagebreak\noindent{\bf Example #1.} }

\define\Rem#1 {\par\medpagebreak\noindent{\bf Remark #1.\ } \ }

%

\define\squash#1#2#3{\vspace{-#1\jot}\intertext{#2}\vspace{-#3\jot}}
\define\eqtext#1#2#3{{\vskip -#1\jot}\noindent#2{\vskip -#3\jot}}

%

\newcount\endnoteno
\newcount \nb
\endnoteno=1
\nb=1
\def\NB{\ifproofmode${}^{\the\nb}$ \global\advance\nb by 1\fi}
\def\myitem{\item{\bf[\the\endnoteno]}\global\advance \endnoteno by 1}

%

\def\qed{\hfill\hbox{\vrule width 4pt height 6pt depth 1.5 pt}}

%
\hyphenation{dif-fer-en-tial di-men-sion-al Helm-holtz
     Cum-mings  Czech-o-sla-vak }

%

\multisectiontrue
\bibmakemodetrue
\pageno=1


\bibmakemodefalse
\bibcitemodetrue
\pagewidth{6.00 truein}
\pageheight{8.25truein}
\vcorrection{.25 truein}

\vcorrection{-.5cm}

\TagsOnRight

\def\sd{\mathbin{ \raise0.0pt\hbox{ \vrule height5pt width.4pt depth0pt}\!\times}}
\def\reals{\text{{\rm I \! \!\! R}}}

\def\Sl{\text{{\bf sl}}}
\def\sol{\text{{\bf s}}}
\def\heis{\text{{\bf n}}}
\def\lag{\text{{\bf g}}}
\def\lo{\text{{\bf o}}}
\def\gl{\text{{\bf gl}}}
\def\bq{\text{{\bf q}}}
\def\bb{\text{{\bf b}}}
\def\ba{\text{{\bf a}}}
\def\be{\text{{\bf e}}}
\def\A5{$\text{{\bf A5}}^*$}
\def\bm{\text{{\bf m}}}
\def\tlah{\tilde{\text{{\bf h}}}}
\def\lah{\text{{\bf h}}}

\def\qed{\hfill\hbox{\vrule width 4pt height 6pt depth 1.5 pt}}

\def\boxit#1{\vbox{\hrule\hbox{\vrule\kern6pt\vbox{\kern4pt#1\kern1pt}\kern1pt }}}

\NoBlackBoxes

\topmatter

\title Non-reductive Homogeneous Pseudo-Riemannian Manifolds of Dimension Four\endtitle

\author M.E. Fels, A.G. Renner  \endauthor
\affil Department of Mathematics and Statistics, Utah State University, Logan Utah, 84322 \endaffil
\address Department of Mathematics and Statistics, Utah State University, Logan Utah, 84322 \endaddress
\email fels\@math.usu.edu,arenner\@cc.usu.edu \endemail
\date October 12, 2004\enddate
\subjclass Primary  53C30 \endsubjclass
\keywords Homogeneous  pseudo-Riemannian, Einstein space \endkeywords
\endtopmatter

\document

\head  Abstract\endhead

\noindent

A method, due to \'Elie Cartan, is used to give an algebraic classification of the non-reductive homogeneous pseudo-Riemannian manifolds of dimension four.
Only one case with Lorentz signature can be Einstein without having constant curvature, and two cases with (2,2) signature are Einstein of which one is Ricci-flat. If a four-dimensional non-reductive homogeneous pseudo-Riemannian manifold is simply connected, then it is shown to be diffeomorphic to $\reals^4$. All metrics for the simply connected non-reductive Einstein spaces are given explicitly. There are no non-reductive pseudo-Riemannian homogeneous spaces of dimension two and none of dimension three with connected isotropy subgroup.

\head 1.  Introduction \endhead

A homogeneous space $G/H$, where $G$ is a Lie group and $H$ a closed Lie subgroup, is {\it reductive }
\cite{kn:1996a} if  the Lie algebra $\lag$ of $G$ may  be decomposed into a vector-space direct sum $ \lag = \lah \oplus \bm$ where $\bm$ is an $Ad(H)$-invariant complement to $\lah$.
If $G/H$ is a reductive homogeneous space which admits a pseudo-Riemannian metric with $G$ acting by isometries, the curvature tensor takes on a particularly simple form. For this reason, the geometry of these spaces has been well studied \cite{kn:1996a} \cite{besse:1987a}, and some classification results are known \cite{gadea:1997a}. On the other hand, little is known about the structure of non-reductive homogeneous pseudo-Riemannian manifolds and the purpose of this paper is classify and investigate the basic geometry and topology of these special manifolds in low dimensions.

While it is easy to construct non-reductive homogeneous spaces, it is quite a bit more difficult to construct examples where $G$ is the isometry group of a pseudo-Riemannian metric on $G/H$. The difficulty is that if $G$ is the isometry group of a Riemannian metric on $G/H$ then $Ad(H)$ is compact, so $G/H$ is automatically reductive (see section 4 for an algebraic proof). Therefore, to construct examples of non-reductive pseudo-Riemannian homogeneous spaces only metrics with indefinite signature need to be considered. These facts are mentioned in \cite{besse:1987a}, but no non-reductive examples are given. In the article \cite{koorn:1981a} the author studies the ring of invariant differential operators on non-reductive homogeneous spaces but only considers geometric examples which turn out to be reductive.

In the book {\it Le\c cons sur la g\'eom\'etrie des espaces de Riemann} \cite{cartan:1951a}, Cartan outlines a method in which questions about the geometry of homogeneous Riemannian manifolds become algebraic questions about Lie algebras. Cartan used his technique to classify the three-dimensional simply connected Riemannian homogeneous spaces which admit a group of isometries of dimension at least 4.
Ishihara \cite{ishihara:1955a} used Cartan's method
to classify the four-dimensional Riemannian manifolds with transitive isometry groups while  Jensen \cite{jensen:1969a} used this technique to determine the simply connected  homogeneous Einstein spaces of dimension 4.
An alternative approach to the classification of low dimensional homogeneous Riemannian
manifolds was given in \cite{bbergery:1978a}, but this approach utilizes the compactness of the isotropy subgroup and so can't be used here.

Cartan's method works perfectly well for pseudo-Riemannian homogeneous spaces. We use this method in section 5 to first show that there are no two or three-dimensional non-reductive homogeneous pseudo-Riemannian manifolds.  We then classify the four-dimensional non-reductive homogeneous pseudo-Riemannian manifolds and show in section 6 that if these four-dimensional homogeneous spaces are simply connected then they are diffeomorphic to $\reals^4$. As a consequence of the calculations in section 5, we identify the cases which are Einstein and compare them with those in \cite{komrakov:2001a}. Finally, in section 7 we construct the corresponding homogeneous Einstein metrics on $\reals^4$ (the simply connected spaces) for the three cases in which they exists.

\initializesection2

\head 2. The Classification and Einstein metrics\endhead

In this section, we provide a summary of the classifications proved
in section 5 and then list the possible Einstein metrics which are
found in section 7 when $G/H$ is assumed to be simply connected.

If $\eta$ is a pseudo-Riemannian metric on the manifold $G/H$ and
$G$ acts effectively and by isometries, we say the pair $(G/H,\eta)$
is a {\it homogeneous pseudo-Riemannian manifold}. We also use the
convention that the bilinear form $\eta$ on an $n$-dimensional
Lorentz manifold has signature $(n-1,1)$. \noindent
\proclaim{Theorem 2.1} Let $(G/H,\eta)$ be a homogeneous Lorentz
manifold. If $G/H$ is two-dimensional, then $G/H$ is reductive. If
$G/H$ is three-dimensional and $H$ is connected, then $G/H$ is
reductive.
\endproclaim
\noindent
Let $\lah $ be a Lie subalgebra of the Lie algebra $\lag$ and denote this pair by $(\lag,\lah)$.

\proclaim{Definition} The Lie algebra pairs $(\lag,\lah)$ and $(\lag',\lah')$ are isomorphic if there exists an isomorphism $\Phi:\lag \to \lag '$ such that $ \Phi(\lah) = \lah
'$.
\endproclaim
\noindent

For every homogeneous space $G/H$, let $\lag$ be the the Lie algebra
of $G$ and $\lah$ the Lie algebra of $H$, let $(\lag,\lah)$ be the
associated Lie algebra pair. In the next theorem, we list all
possible non-isomorphic Lie algebra pairs for the non-reductive
four-dimensional homogeneous spaces that are classified in section
5. We use the table of Lie algebras in \cite{psw:1976a} and refer to
these algebras by $A_{x,y}$, as given on page 990.

\noindent
\proclaim{Theorem 2.2} Let $(G/H,\eta)$ be a four-dimensional homogeneous Lorentz manifold where $H$ is connected. If $G/H$ is not reductive, then the Lie algebra pair $(\lag ,\lah)$ is isomorphic to one in the following list.
\roster

\smallskip
\item"{\bf A1}"  The Lie algebra $\lag$ is the decomposable five-dimensional algebra $\Sl(2,\reals) \oplus \sol(2)$, where $\sol(2)$ is the two-dimensional solvable algebra. There exists a basis for $\lag$ where the non-zero products are
$$
[\be_1, \be_2]=2 \be_2,\ [\be_1,\be_3]=-2\be_3,\ [\be_2,\be_3]=\be_1,\ [\be_4
,\be_5]= \be_4 .
$$
The isotropy is $\lah = span \{ \be_3 + \be _4 \}$.

\smallskip
\item"{\bf A2}"  The Lie algebra $\lag$ is the one-parameter family of five-dimensional solvable Lie algebras $ A_{5,30}$. There exists a basis for $\lag$ where the non-zero products are
$$
\alignat 3
[\be_1,\be_5]  &= (\alpha+1)\be_1, \quad
&[\be_2,\be_4] &= \be_1 , \quad
&[\be_2,\be_5] &=\alpha\be_2 , \\
[ \be_3,\be_4] &=\be_2,\quad
&[\be_3,\be_5] &=(\alpha -1) \be_3,\quad
&[\be_4,\be_5] &=\be_4 ,
\endalignat
$$
where all values of $\alpha \in \reals$ are admissible. The isotropy is $\lah =  span\{ \be_4 \}$.

\smallskip
\item"{\bf A3}" The Lie algebra $\lag$ is one of the five-dimensional solvable algebras $A_{5, 37}$ or $A_{5,36}$. There exits a basis for $\lag$ where the non-zero products are
$$
[\be_1,\be_4]=2 \be_1\ ,\ [\be _2,\be_3]=\be_1,
[\be_2,\be_4]=\be_2,\
[\be_2,\be_5]=-\epsilon \be_3, \
[\be_3,\be_4]=\be_3, \
[\be_3,\be_5]=\be_2 \, ,
$$
where $\epsilon = 1 $ for $A_{5,37}$ and $\epsilon = -1$ for $A_{5,36}$. The isotropy is $\lah = span \{ \be_3 \} $.

\smallskip
\noindent
\item"{\bf A4}" The Lie algebra $\lag$ is the six-dimensional algebra $\Sl(2,\reals) \sd  \heis (3) $ where $ \heis(3)$ is the three-dimensional Heisenberg algebra. There exits a basis for $\lag$ where the non-zero products are
$$
\alignat 4
 [\be_1, \be_2]  & =2 \be_2,\quad & [\be_1,\be_3] &  =-2 \be_3,\quad & [\be_2,\be_3] &= \be_1, \quad & [\be_1, \be_4] & = \be_4,  \\
 [\be_1 ,\be_5]& =-\be_5, &[\be_2,\be_5]& = \be_4, & [\be_3, \be_4]&= \be_5, &[\be_4,\be_5]&=\be_6.
\endalignat
$$
The isotropy is ${\lah} = span \{ \be_3 + \be_6\, \be_5\}$. The algebra is sometimes called the Schroedinger algebra.

\smallskip
\noindent
\item"{\A5}" The Lie algebra $\lag$ is the seven-dimensional algebra $\Sl (2,\reals) \sd A_{4,9}^1 $. There exists a basis for $\lag$ where the non-zero products  are
$$
\alignat 5
[\be_1, \be_2]  & =2 \be_2,\quad & [\be_1,\be_3] &  =-2 \be_3,\quad & [\be_1,\be_5] &=- \be_5, \quad & [\be_1, \be_6] & = \be_6,\quad & [\be_2 ,\be_3]& =\be_1,  \\
[\be_2,\be_5]& = \be_6, & [\be_3, \be_6]&= \be_5, &[\be_4,\be_7]&= 2 \be_4, & [\be_5,\be_6]& =\be_4,& [\be_5,\be_7]&=\be_5 \,\quad [\be_6,\be_7]=\be_6 .
\endalignat
$$
The isotropy is ${\lah} = span \{ \be_1 + \be_7,\be_3- \be_4,\be_5 \}$.
\endroster
\endproclaim

In the book Einstein Spaces \cite{petrov:1969a}, Petrov gave a
fairly comprehensive list of the possible infinitesimal generators
for the isometry algebras of a four-dimensional Lorentz manifold.
The Lie algebras in {\bf A1} and {\bf A4} should appear on the list,
but they don't.

We now list the possibilities when the signature is $(2,2)$.

\noindent
\proclaim{Theorem 2.3} Let $(G/H,\eta)$ be a four-dimensional homogeneous pseudo-Riemannian manifold of signature (2,2) where $H$ is connected. If $G/H$ is not reductive, then the Lie algebra pair $(\lag ,\lah)$ is isomorphic to one in the following list.
\roster

\smallskip
\item"{\bf A1}"\!\!\! {\bf -A3}  The corresponding Lie algebra pairs in Theorem 2.2.

\smallskip
\item"{\bf B1}" The Lie algebra $\lag$ is the five-dimensional algebra $ \Sl(2,\reals)\sd \reals^ 2$. There exists a basis for $\lag$ where the non-zero products are
$$
\alignat 4
 [\be_1, \be_2]  & =2 \be_2,\quad & [\be_1,\be_3] &  =-2 \be_2,\quad & [\be_2,\be_3] &= \be_1, \quad & [\be_1, \be_4] & = \be_4,  \\
 [\be_1 ,\be_5]& =-\be_5, &[\be_2,\be_5]& = \be_4, & [\be_3, \be_4]&= \be_5. &&  \
\endalignat
$$
The isotropy is $\lah = span \{ \be _3\}$.

\smallskip
\item"{\bf B2}" The Lie algebra $\lag$ is the six-dimensional Schroedinger algebra $\Sl(2,\reals) \sd  \heis (3) $ in {\bf A4} of Theorem 2.2 and the isotropy ${\lah} = span \{ \be_3 -\be_6 , \be_5\}$.

\smallskip
\item"{\bf B3}" The Lie algebra $\lag$ is the six-dimensional algebra $ \Sl(2,\reals) \sd  \reals^2 \oplus \reals $. There exits a basis for $\lag$ where the non-zero products are
$$
\alignat 4
 [\be_1, \be_2]  & =2 \be_2,\quad & [\be_1,\be_3] &  =-2 \be_3,\quad & [\be_2,\be_3] &= \be_1, \quad & [\be_1, \be_4] & = \be_4,  \\
[\be_1 ,\be_5]& =-\be_5, &[\be_2,\be_5]& = \be_4, & [\be_3, \be_4]&= \be_5.
&&
\endalignat
$$
The isotropy is ${\lah} = span \{\be_3, \be_5 + \be_6\}$.

\endroster
\endproclaim

The following theorem proved in section 6, gives a complete classification when the space is simply connected.

\proclaim{Theorem 2.4} Let $G/H$ be a simply connected non-reductive pseudo-Riemannian homogeneous space of dimension four, then
\roster

\item"{\rm i]}" $G/H$ is diffeomorphic to $\reals^4$, and

\item"{\rm ii]}" if $G$ is the full isometry group then the Lie algebra pair for $G/H$ is equivalent to one in Theorem 2.2 excluding \A5, or to one in Theorem 2.3.

\endroster
Conversely, for any Lie algebra pair from Theorem 2.2 except \A5, or any in Theorem 2.3, there exists a pseudo-Riemannian metric on $\reals^4$ (subject to the conditions on the signature), where the isometry group acts transitively on $\reals^4$, the Lie algebra of the isometry group is the given Lie algebra $\lag$, and the Lie algebra of the isotropy at a point is (conjugate to) $\lah$.
\endproclaim

We show in Lemma 5.1 that only {\bf A2} in Theorem 2.2 or 2.3 with
$\alpha=\frac{2}{3}$ and {\bf B3} lead to Einstein spaces which are
not of constant curvature. Furthermore, {\bf B3} is Ricci-flat. By
using this result we prove in section 7  the following theorem which
gives complete list of all the homogeneous Einstein metrics which
are not of constant curvature for the simply connected non-reductive
pseudo-Riemannian manifolds of dimension 4.

\proclaim{Theorem 2.5} Let $(G/H,\eta)$ be a simply connected
non-reductive homogeneous space of dimension 4 which is Einstein and
not of constant curvature. If $\eta$ is Ricci-flat, then the Lie
algebra pair is isomorphic to {\bf B3} and there exists global
coordinates $(y^i)_{i=1..4}$ on $G/H=\reals^4$ such that the metric
is
$$
\eta = 2 e^{y^1}\cos y^2( dy^1 d y^4 -dy^2  dy^3) -2 e^{y^1}\sin y^2
(dy^1  dy^3 + dy^2 dy^4)+L e^{4y^1} (dy^2)^2
$$
for some $L \in \reals^*$. Otherwise the Lie algebra pair is
isomorphic to {\bf A2} with $\alpha =2/3$ and there exists global
coordinates $(y^i)_{i=1..4}$ on $G/H=\reals^4$ such that the metric
is
$$
\eta = a_1 e^{ -\textstyle{\frac{4}{3}} y^4}(2 dy^1  dy ^3 -
(dy^2)^2 ) + a_2 e^{\frac{2}{3}y^4} (dy^3)^2 +2 a_3
e^{\frac{1}{3}y^4} dy^3 dy^4 + a_4 (dy^4)^2
$$
for some choice of $a_i \in \reals, \ i=1..4 $ where $a_1 a_4 \neq
0$, and $a_2\neq 0$.
\endproclaim
It is worth noting that determining the Lie algebra of the isometry
group for the Ricci flat metrics in this theorem is non-trivial.

\initializesection3

\head 3.  Cartan's Approach to the Geometry of Homogeneous Spaces \endhead

Let $\eta^0$  be a non-degenerate symmetric bilinear form on
$\reals^n$ with signature $(p,\tilde p)$, and $O(p,\tilde p) \subset
GL(n,\reals)$ be the Lie group preserving $\eta^0$. Let $ (M,\eta) $
be a pseudo-Riemannian manifold of signature $(p,\tilde p)$, and
$\pi:O(M) \to M$ be the orthonormal frame bundle corresponding to
$\eta^0$ defined by
$$
O(M) = \{ u_p:\reals^n \to T_p M \ | \quad \eta(u_p(X),u_p(Y)) = \eta^0(X,Y) \} \ .
$$
Denote the right action of $a\in O(p,\tilde p)$ on $u\in O(M)$ by $ u \, a $, and for $X \in \lo (p,\tilde p)$,  let $\tilde X$ be the corresponding infinitesimal generator on $O(M)$ defined by
$$
{\tilde X}_u =  \frac{d\ }{dt} \left( u \, exp(tX_e) \right) |_{t=0}
\ .
$$
The canonical  $\reals^n$-valued one-form $\theta$ and the $\lo(p,\tilde p)$-valued
connection one-form $\omega$ on $O(M)$ have the following properties \cite{kn:1996a} (pp. 118-121)
$$
\iota_Z \theta  = u^{-1} \pi_* (Z), \quad
 \iota_{\tilde X} \omega  = X, \quad
 d \theta = - \omega\wedge \theta,
\Tag{21a}
$$
\noindent where $X \in \lo(p,\tilde p)$, $Z \in T_u(O(M))$ and
$\iota$ is left interior multiplication. The $\lo(p,\tilde
p)$-valued curvature two-form $ \Omega = d \omega + \omega \wedge $
satisfies
$$
\iota _{\tilde X}\Omega =0 ,\quad    \Omega \wedge \theta = 0  , \quad d\Omega = \Omega \wedge \omega - \omega \wedge \Omega\,  .
\Tag{21d}
$$
The forms $\theta, \ \omega,$ and $ \Omega$ satisfy the equivariance conditions
$$
a^* \theta = a^{-1}(\theta) , \quad   a^*\omega = Ad(a^{-1}) \omega, \quad  a^*\Omega = Ad(a^{-1}) \Omega \qquad \text{where} \ a \in O(p,\tilde p) .
$$
If $H$ is connected, then \eq{21a} and \eq{21d} imply the
equivariance of $\theta$, $\omega$ and $\Omega$.
\smallskip

Let $g:M\to M$ be an isometry of the pseudo-Riemannian manifold $(M,\eta)$, and $\phi_g$ be the lift of the diffeomorphism $g$ of $M$ to the frame bundle,
$$
\phi_g(u) = g_* u \ , \qquad u \in F(M).
$$
Since $g$ is an isometry, the subset $O(M) \subset F(M)$ is
invariant under $\phi_g$. The forms $\theta,\omega,$ and $\Omega$
satisfy the invariance properties
$$
\phi_g^* \theta = \theta \ , \ \phi_g^*\omega =  \omega \ , \ \phi_g^*\Omega = \Omega \  .
\Tag31
$$

Suppose that $(G/H,\eta)$ is a homogeneous pseudo-Riemannian
manifold and let $\sigma =[H] \in G/H$, and $u_\sigma \in O(G/H)$ be
an orthonormal frame at $\sigma$. The linear isotropy representation
$ \rho:H \to O(p,\tilde p)$ is defined by
$$
u_{\sigma} \rho(h) = (L_h)_* u_{\sigma}
\Tag32
$$
where $L_h$ is right multiplication in $G$ by $h \in H$. The
differential of $\rho$ defines a homomorphism $ \rho_* :\lah \to
\lo(p,\tilde p) $. Since $G$ acts effectively and by isometries, the
linear isotropy representation of $H$ is faithful. Following Cartan
\cite{cartan:1951a} (or see Jensen \cite{jensen:1969a}), we define
the function $\Psi:G \to O(G/H)$ by
$$
\Psi(g) = g_* u_{\sigma}
\Tag33
$$
which makes the diagram
$$
\CD
G        @> \Psi >> O(G/H) \\
{\bq} @VVV        @VVV \!\!\!\!\!\!\! \pi \\
G/H        @=  G/H \ .
\endCD
\Tag34
$$
commutative.
The map $\Psi$ is equivariant with respect to the left action of $G$ on $G$ and the
action of $G$ on $O(M)$.
It is also equivariant with respect to the linear isotropy representation. Therefore, $\Psi$
satisfies
$$
\Psi( g h )  =  \Psi(g) \rho(h) ,  \quad \text{\rm and} \quad
\Psi(g_1 g_2 )  = \phi_{g_1} \circ \Psi(g_2) \ . \Tag35
$$
\noindent
By defining the forms
$$
\hat \theta = \Psi^* \theta, \quad\hat \omega = \Psi^* \omega,\quad
\hat \Omega=\psi^* \Omega
$$
which are $G$-invariant on account of the equivariance  of $\Psi$
and  equation \equationlabel{3}{31}{}, we obtain the following
structure on the Lie algebra $\lag$ of $G$.

\proclaim{Lemma 3.1} Let $(G/H,\eta)$ be an n-dimensional homogeneous pseudo-Riemannian manifold with Lie algebra pair $(\lag,\lah)$.
There exists an injective homomorphism $\rho_* : \lah \to \lo(p,\tilde p)$, an $\reals ^n$-valued one-form $\hat \theta : \lag \to \reals^n$, and an $\lo(p,\tilde p)$-valued one-form $\hat \omega : \lag \to \lo(p,\tilde p)$ satisfying
\smallskip
$1]$ $ \ker\hat \theta= \lah$,

\smallskip
$2]$ $ \hat \omega(X) = \rho_*(X) $, and

\smallskip
$3]$ $ d \hat \theta = -\hat \omega  \wedge  \hat\theta \ $

\smallskip
\noindent
where $X \in \lah$.
Furthermore, the $\lo(p,\tilde p)$-valued two-form
$ \Omega = d \hat \omega + \hat \omega \wedge \hat \omega $ satisfies

\smallskip
$4]$ $\iota _X \Omega = 0 $, \quad $\hat \Omega \wedge \hat \theta = 0 $, and \quad $d \hat\Omega =  \hat \Omega \wedge \hat\omega -\hat \omega \wedge \hat \Omega $.
\endproclaim


\noindent Lemma 3.1 has the following partial converse.

\proclaim{Lemma 3.2} Let $\lah$ be a Lie algebra, $\rho_*:\lah \to
\lo(p,\tilde p)$ a monomorphism, and $\lag$ be the vector space
$\reals^n \oplus \lah$. Suppose there exists forms $\theta : \lag
\to \reals ^n $, $\omega : \lag \to \lo(p,\tilde p) $, and $\Omega:
\lag \wedge \lag \to \lo(p,\tilde p)$ satisfying

\medskip
$\hat 1]$ $ \ker \theta = \lah $,

\medskip
$\hat 2]$ $ \omega(X) = \rho_*(X) $,

\medskip
$\hat 3]$ $\iota_X \Omega = 0 $, \ and  $  \Omega \wedge \theta = 0 $.

\medskip
\noindent
If we define $d \theta = - \omega \wedge \theta $ and $\Omega$ satisfies

\medskip
$\hat 4]$ $d \Omega =  \Omega \wedge \omega -\omega \wedge \Omega $,

\medskip
\noindent then with $d\omega = \Omega - \omega \wedge \omega$,
$\lag$ is a Lie algebra where $\alpha([X,Y])=-d\alpha(X,Y), \alpha
\in \lag^*, X,Y \in \lag$.
\endproclaim

The principle step in Cartan's approach to the classification of homogeneous pseudo-Riemannian manifolds is to start with a subalgebra $\lah \subset \lo(p,\tilde p)$ and then classify all Lie algebras that satisfy Lemma 3.2.  To simplify this classification, one expects that we only need the conjugacy class of the subalgebra $\lah \subset \lo(p,\tilde p)$ under inner automorphism, but
slightly more is true.

\proclaim{Lemma 3.3} Let $ \lah$ and $\tlah $ be two Lie algebras, and let $\rho_*:\lah \to \lo(p,\tilde p)$ and $\tilde \rho_* : \tlah \to \lo(p,\tilde p)$ be monomorphisms. Suppose there exists an inner automorphism $\psi:\gl(n,\reals) \to \gl(n,\reals)$ which restricts to an isomorphism $\phi:\lah \to \tlah$ such that
$$
\tilde \rho_* (\phi(X) ) = \psi (\rho_*(X)) \quad X\in \lah \ .
$$
Then the pairs $(\lag,\lah)$ which satisfy Lemma 3.2 are in one to one correspondence with the pairs $(\tilde {\lag}, \tilde { \lah })$ which satisfy Lemma 3.2.
\endproclaim

\demo{Proof} Suppose the inner automorphism $\psi$ is conjugation by
the matrix $A\in GL(n,\reals)$. It is then easy to check that the
vector-space isomorphism $T:\reals^n \oplus \lah \to \reals \oplus
\tlah$  defined by
$$
T(\xi,X) = (A\xi, \phi(X))
$$
provides a correspondence. \qed
\enddemo


\noindent
Note that every inner automorphism of $\lo(p,\tilde p)$ satisfies this lemma.

\initializesection4

\head 4. Non-Reductive Homogeneous Spaces \endhead

The preceding section described Cartan's procedure for constructing all
possible isomorphism classes of Lie algebra pairs $(\lag, \lah)$ for
homogeneous pseudo-Riemannian manifolds by starting from the
inequivalent subalgebras of $\lo(p,\tilde p)$ under the automorphisms described in Lemma 3.3.
In principle, one could find a general classification of the four-dimensional homogeneous pseudo-Riemannian manifolds starting with the entire list of subalgebras for the Lie algebras $\lo(3,1)$ and $\lo(2,2)$. This classification would be rather daunting
because the known lists of inequivalent subalgebras under inner automorphisms are quite large \cite{pwz:1975a}, \cite{psw:1977a}.
In this section, we simplify the classification problem by proving a lemma which reduces the possible subalgebras $\lah \subset \lo(p,\tilde p)$ associated with a non-reductive homogeneous pseudo-Riemannian manifold $(G/H,\eta)$.

We start with the following characterization of reductive
homogeneous spaces see \cite{kn:1996a}, p. 103, Theorem 11.1.

\proclaim{Lemma 4.1} A homogeneous space $G \to G/H$ is reductive if and only if the principal $H-$bundle $G \to G/H$ admits a $G$-invariant connection.
\endproclaim
\noindent

The following lemma greatly simplifies the classification problem.

\proclaim{Lemma 4.2} If $G/H$ is a pseudo-Riemannian homogeneous space and
 $O(p,\tilde p)/\rho(H) $ is a reductive homogeneous space, then $G/H$ is reductive.
\endproclaim

\demo{Proof} Let $\tilde {\lah}$ be the Lie algebra of $\rho(H)$, and  $\lo(p,\tilde p) =\tilde {\lah} \oplus \bm $ be a reductive decomposition of $\lo(p,\tilde p)$. Decompose the connection form on $O(G/H)$ as $\omega = \omega_{\tilde {\lah}} + \omega_{\bm} $ where $\omega_{\tilde {\lah}} $ takes values in $\tilde {\lah}$ and $\omega_{\bm}$ takes values in $\bm$.  By using the map $\rho$ defined in \equationlabel{3}{32}{} and $\Psi $ defined
in \equationlabel{3}{33}{}
we prove that the $\lah$-valued form $\rho^{-1}_* \circ (\Psi^* \omega_{\tilde {\lah}})$ defines a $G$-invariant connection on $G/H$.

The $G$-invariance of $\rho^{-1}_* \circ (\Psi^* \omega_{\tilde
{\lah}})$ follows from the equivariance of $\Psi$ in
\equationlabel{3}{35}{} together with \equationlabel{3}{31}{}. In
order that this form defines a connection we need to check that the
two conditions on p.64 in \cite{kn:1996a} are satisfied. To check
the first condition on p.64 \cite{kn:1996a},  we use
\equationlabel{3}{21a}{} and compute
$$
\rho^{-1}_* \circ (\Psi^* \omega_{\tilde {\lah}} (X_e)) = \rho^{-1}_* \circ \omega_{\tilde {\lah}} ( \rho_*(X_e) ) = X_e \ .
$$
This verifies condition one.
We now check the second condition. It follows from the hypothesis in the lemma and the equivariance of the connection form $\omega$ that
$$
R^*_a \omega_{\tlah } = Ad_{a^{-1}} \omega_{\tlah } \quad \text{\rm and}  \quad
R^*_a \omega_{\bm } = Ad_{a^{-1}} \omega_{\bm} .
$$
Now from the $H$-equivariance of $\Psi$ in \equationlabel{3}{35}{}, the equation above, and the identity $\rho^{-1}_* \circ Ad _{\rho(h)} = Ad_h \circ \rho^{-1}_* $ it follows that
$$
R_h^* \left( \rho^{-1}_* \circ (\Psi^* \omega_{ \tilde{\lah}}) \right) =
 \rho^{-1}_* \circ \left(\Psi^*  \rho(h)^*\omega_{\tilde{ \lah}} \right) =
 \rho^{-1}_* \circ \left(Ad _{\rho(h)^{-1}} \circ \Psi^* \omega_{\tilde {\lah}} \right) =
 Ad_{h^{-1}} \circ \rho^{-1}_* \circ \left(\Psi^* \omega_{\tilde {\lah}} \right).
$$
This verifies the second condition. Therefore, the $\lah$-valued
form $\rho^{-1}_* \circ \Psi^* \omega_{\tlah}$, is a $G$-invariant
connection on $G/H$ and by Lemma 4.1, $G/H$ is reductive. \qed
\enddemo
The proof of Lemma 4.2 is similar to the proof of Proposition 6.4 on
pg. 83 in \cite{kn:1996a}. This proposition implies that if
$O(p,\tilde p)/H$ is reductive, then the metric connection is
reducible to $H$.

This lemma has a few simple but interesting corollaries.

\proclaim{Corollary 4.1} If $\rho_*(\lah) \subset \lo(p,\tilde p)$ is a non-degenerate subspace with respect to the Killing form of $\lo(p,\tilde p)$, then $G/H$ is reductive.
\endproclaim

\noindent
Lemma 4.2 also provides an algebraic proof of the following corollaries.

\proclaim{Corollary 4.2} If $G/H$ admits a $G$ invariant Riemannian metric, then $G/H$ is a reductive homogeneous space.
\endproclaim

\proclaim{Corollary 4.3} If $(G/H,\eta)$ is a two dimension homogeneous Lorentz manifold, then $G/H$ is reductive.
\endproclaim

\initializesection5
\head 5. The Computations .\endhead

By starting with the inequivalent subalgebras of $\lo(2,1)$, $\lo(3,1) $, and $\lo(2,2)$ we prove  Theorems 2.1, 2.2 and 2.3 by building
all non-reductive Lie algebra pairs $(\lag, \lah)$ which satisfy Lemma 3.2.
All inequivalent subalgebras of $\lo(2,1)$, $\lo(3,1) $, $\lo(2,2)$ under inner automorphisms are known. Although this list is rather long, Lemma 4.2 says that we need only those subalgebras that are not reductive in their respective algebras. With this reduced list of subalgebras, the equivalence problem in Lemma 3.3 is much easier.
By using this final list of inequivalent subalgebras, we determine those which extend to a Lie algebra that satisfies Lemma 3.2 and not Lemma 4.1 (see also Lemma A1). The resulting Lie algebra pairs are then put into a canonical form which proves Theorems 2.1, 2.2 and 2.3.
\smallskip

Let $ (e^i_{j}) $ denote the standard basis for $gl(n,\reals)$  where
$$
(e^i_{j}) ^k_l = \delta^{i}_{k}\delta^l_j \ ,\quad 1 \leq i,j,k,l \leq n .
$$
Hereafter we omit writing the isomorphism $\rho_*$ between $\lah$ with basis $\{ \be_\alpha \}_{\alpha=1\ldots q}$ and $\rho_*(\lah) \subset \lo(p,\tilde p)$ with basis $\{ \bb_\alpha \}_{\alpha=1\ldots q}$. Given two differential one forms $\sigma^1 , \sigma ^2 \in \Omega^1(M)$, we use
the convention
$$
\sigma^1 \sigma^2 = \frac{1}{2}( \sigma^1 \otimes \sigma^2 + \sigma^2 \otimes \sigma^1)
$$
for the symmetric tensor product. Other notation that is used is given in Appendix A.

\smallskip
\noindent
$\underline{\text{\it Proof of Theorem 2.1 - $\lo(2,1)$:}}$ Let $\{ \sigma^i \}_{i=1\ldots 3 }$ denote the standard basis for $(\reals^3 )^*$, and
$$
\eta^0 = (\sigma^1)^2 + 2 \sigma^2  \sigma ^3 .
$$
For $\lo(2,1)$ we use the basis
$$
B_1 = e^2_2-e^3_3, \ B_2 = e^1_2 - e^3_1, \ B_3= e^2_1- e^1_3 \ .
$$

\smallskip
\noindent
Of the inequivalent subalgebras of $\lo(2,1)$ under inner automorphism, only two are not reductive.  In each case, by using equations (A.1) and (A.3), equation (A.5) always has a solution, so for these two subalgebras, the constructed homogeneous space will be reductive. Here are the details.

\smallskip
\noindent ${\text{\bf Case 1:}}$ The isotropy subalgebra is $\lah =
span\{\bb_1 = B_3\} $. By using the basis $\{ \tilde {\bb}_1 = B_1 \
, \tilde {\bb}_2 = B_2, \bb_1 \}$  for $\lo(2,1)$, equations (A.1)
and (A.3) give
$$
\tilde \omega^1 =  p_1 \theta^1 +p_2 \theta^3, \quad \tilde \omega^2= -p_1 \theta^3\ .
$$
Equation (A.5) has the general solution
$$
r^1_1 =-p_2, \quad r^1_2=p_1,
$$
and $(\omega ^ 1 -p_2 \theta^1 +p_1 \theta^2)\otimes \be_1$ defines a $G$-invariant connection.
\smallskip
\noindent
${\text{\bf Case 2:}}$ The isotropy subalgebra is $\lah = span\{\bb_1 = B_1, \bb_2 = B_3\} $. By using the basis  $\{\tilde {\bb}_1 = B_2,  \bb_1, \bb_2 \}$  for $\lo(2,1)$, equations (A.1) and (A.3) give
$$
\tilde \omega^1 =  p_1 \theta^3.
$$
Equation (A.5) has the general solution
$$
r^1_1 =p_1, \quad r^1_2=r^1_3=r^2_1=r^2_3=0, \quad r^2_2 = -p_1,
$$
and $(\omega ^ 1 +p_1 \theta^1)\otimes \be_1 + (\omega ^ 2 -p_1 \theta^2)\otimes \be_2 $ defines a $G$-invariant connection.
\smallskip
\noindent
This proves Theorem 2.1. \qed

\medskip
\noindent
$\underline{\text{\it Proof of Theorem 2.2 - $\lo(3,1)$:}}$  Let $\{ \sigma^i \}_{i=1\ldots 4}$ denote the standard basis for $(\reals^4)^*$, and
$$
\eta^0 = (\sigma^1)^2 +(\sigma^2)^2 + 2 \sigma^3  \sigma ^4 .
$$
The basis for $\lo(3,1)$ we use is
$$
B_1 = e^2_1 - e^1_2, \quad
B_2 = e^4_4 - e^3_3, \quad
B_3 = e^4_1 - e^1_3, \quad
B_4 = e^4_2-e^2_3, \quad
B_5 = e^1_4-e^3_1, \quad
B_6= e^3_2 -e^2_4. \quad
$$
\smallskip
\noindent
On page 1605 of \cite{pwz:1975a}, the inequivalent subalgebras of $\lo(3,1)$, under inner automorphisms are listed. Of these subalgebras,
labeled $F_1$ to $F_{15}$, seven are not reductive in $\lo(3,1)$.

\medskip
\noindent
${\text{\bf Case 1:}}$ We consider the non-reductive subalgebras  of $\lo(3,1)$ which admit a solution to equation (A.5). Therefore they always lead to a reductive homogeneous space.

\smallskip
\noindent
${\text{\bf Subcase 1.1:}}$ The subalgebra $F_2$ in \cite{pwz:1975a} is $\lah = span \{\bb_1= B_1,\bb_2= B_2,\bb_3= B_3,\bb_4=B_4\}$.
By using the basis $\{\tilde {\bb}_1=B_5,\tilde{ \bb}_2=B_6, \bb_1,\bb_2,\bb_3,\bb_4\}$  for $\lo(3,1)$, equations (A.1) and (A.3) give
$$
\tilde \omega^1 = 0, \quad \tilde \omega^2=0.
$$
Equation (A.5) has the general solution $\{r^\alpha_i = 0 \}_{\alpha = 1\ldots 2, i= 1\ldots 4}$, and $\omega^\alpha \otimes \be_\alpha$ defines a $G$-invariant connection.

\smallskip
\noindent
${\text{\bf Subcase 1.2:}}$ The subalgebra $F_5$ in \cite{pwz:1975a} is $\lah = span \{\bb_1= \cos \theta \, B_1 +\sin \theta B_2, \bb_2= B_3,\bb_3= B_4\}$ where $ \theta\in (0,\pi), \theta\neq \pi/2 $. By using the basis $\{\tilde {\bb}_1=-\sin \theta B_1 + \cos \theta B_2 ,\tilde{ \bb}_2=B_5, \tilde {\bb}_3=B_6, \bb_1, \bb_2, \bb_3\}$  for $\lo(3,1)$, equations (A.1) and (A.3) give
$$
\tilde \omega^1 =  0, \quad \tilde \omega^2=0 ,\quad \tilde \omega^3 = 0.
$$
Equation (A.5) has the general solution $\{ r^\alpha_i = 0\}_{\alpha=1\ldots 2,i = 1\ldots 4}$ and $ \omega^\alpha \otimes \be_\alpha$ defines a $G$-invariant connection.

\smallskip
\noindent
${\text{\bf Subcase 1.3:}}$ The subalgebra $F_6$ in \cite{pwz:1975a} is
$\lah = span \{\bb_1= B_1, \bb_2= B_3, \bb_3= B_4\}$.
By using the basis $\{\tilde {\bb}_1=B_2,\tilde{ \bb}_2=B_5, \tilde {\bb}_3=B_6,\bb_1,\bb_2,\bb_3\}$  for $\lo(3,1)$, equations (A.1) and (A.3) give
$$
\tilde \omega^1 =  p_1 \theta^3, \quad \tilde \omega^2=0\ ,\quad \tilde \omega^3 = 0.
$$
Equation (A.5) has the general solution
$$
r^1_3 =r^3_1 =t, \quad
r^2_2 =-t, \quad
r^1_1 =r^1_2 =r^1_4 =r^2_3 =r^2_4 =r^3_3 =r^3_4 =0,\quad
r^2_1 =r^3_2 =-p_1 \quad
$$
where $t \in \reals$, and $ (\omega^1+t\, \theta^3) \otimes \be_1 + (\omega^2 -p_1 \theta^1-t\, \theta^2) \otimes \be_2 + (\omega^3 +t\, \theta^1 -p_1 \theta^2 ) \otimes \be_3 $
defines a $G$-invariant connection for any choice of $t\in \reals$.
\smallskip
\noindent

\smallskip
\noindent
${\text{\bf Subcase 1.4:}}$ The subalgebra $F_8$ in \cite{pwz:1975a} is  $\lah = span \{\bb_1 = B_2,\bb_2 = B_3\}$.  By using the basis $\{ \tilde{\bb}_1=B_1,\tilde {\bb}_2=B_4,\tilde {\bb}_3=B_5,\tilde {\bb}_4=B_6,\bb_1,\bb_2 \}$  for $\lo(3,1)$, equations (A.1) and (A.3) give
$$
\tilde \omega^1  = - p_1 \theta ^1, \quad
\tilde \omega^2  = p_1 \theta^4, \quad
\tilde \omega^3  = p_2 \theta^3, \quad
\tilde \omega^4  = p_1 \theta ^3.
$$
Equation (A.5) has the general solution
$$
r^1_1 =p_2,\quad
r^2_1 =r^1_2 =r^2_2 =r^1_3 =r^2_3 =r^1_4 =0, \quad
r^2_4 =-p_2,
$$
and $(\omega^1+p_2\theta^1) \otimes \be_1 + (\omega^2 -p_2 \theta^4) \otimes \be_2 $ defines a $G$-invariant connection.

\medskip
\noindent
We now consider the cases where condition (A.5) is not automatically satisfied.
\smallskip
\noindent
${\text{\bf Case 2:}}$ The subalgebra $F_7$ in \cite{pwz:1975a} is $\lah = span \{\bb_1 = B_2,\bb_2 = B_3,\bb_3=B_4\} $.
By using the basis $\{\tilde{\bb}_1=B_1,\tilde{\bb}_2=B_5,\tilde{\bb}_3=B_6,\bb_1,\bb_2,\bb_3\}$  for $\lo(3,1)$, equations (A.1) and (A.3) give
$$
\tilde \omega^1 = p_1 \theta ^1 +p_2 \theta^2, \quad
\tilde \omega^2 =-p_2 \theta^3, \quad
\tilde \omega^3 =-p_1 \theta^3.
$$
From condition (A.5), $G/H$ is reductive if and only if $ p_1 = p_2 = 0 $, so we assume this is not satisfied. Let $ K = {{p_1}^2+{p_2}^2} $ (which is non-zero). The Bianchi identities give
$$
C^1_{42}=-K,\ C^2_{43}=-K,\ C^3_{41} =-K\, ,
$$
and $ C^\alpha_{jk} =0 $ otherwise. The curvature form is $ \Omega_{ij} = -K \theta_i \wedge \theta_j$, and the homogeneous space will be of constant curvature. The change of basis
$$
\alignat 4
\alpha^1 &= (p_1 \theta^2-p_2\theta^1-\omega^1)/2, \ &
\alpha^3 &=\theta^4 +K^{-1}(p_1 \omega^3-p_2 \omega^2), \ &
\alpha^4 &=\theta^4 +K^{-1}(p_2\omega^2-p_1\omega^3),   &
\alpha^2 &= \! K \theta^3, \\
\alpha^5 &= \sqrt{2}K^{-1}(p_1 \omega^2 + p_2 \omega^3),\ &
\alpha^6 &=\sqrt{2}( p_1 \theta^1+p_2 \theta^2), \ &
\alpha^7 &=(p_2 \theta^1-p_1\theta^2-\omega^1)/2
\endalignat
$$
for $\lag^*$ leads to the multiplication table {\A5} in Theorem 2.2 with isotropy in the dual basis $\lah= span\{ \be_1 +\be_7,\be_3-\be_4,\be_5 \}$.

\smallskip
\noindent
${\text{\bf Case 3:}}$ The subalgebra $F_{10}$ in \cite{pwz:1975a} is $\lah = span \{\bb_1 = B_3,\bb_2 = B_4\} $. By using the basis $\{\tilde{\bb}_1=B_1, \tilde{\bb}_2=B_2, \tilde{\bb}_3=B_5, \tilde{\bb}_4=B_6, \bb_1, \bb_2\}$  for $\lo(3,1)$, equations (A.1) and (A.3) give
$$
\tilde \omega^1 = p_1 \theta^1 + p_2 \theta^2 + p_3 \theta^3, \quad
\tilde \omega^2 = p_2 \theta^1 - p_1 \theta^2 + p_4 \theta^3, \quad
\tilde \omega^3 =-p_2 \theta^3, \quad
\tilde \omega^4 =-p_1 \theta^3.
$$
From condition (A.5), $G/H$ is reductive if and only if $ p_1 = p_2 = 0 $. Let $ K =-( {{p_1}^2+{p_2}^2})$ (which is non-zero) and $C_{13}^1=L(p_1{}^2+4p_2{}^2)$. The Bianchi identities give
$$
\alignat 7
C_{12}^1 &= -3 p_2p_3, &\
C_{14}^1 &= -K, &\
C_{34}^1 &= 2p_4p_2-p_1p_3, &\
C_{24}^1 &= C_{14}^2= 0, &\
C_{23}^1 &= -3Lp_2p_1, &\
p_3 &= 0,
\\
C_{12}^2 &=3 p_1p_3, \ &
C_{24}^2 &= -K, \ &
C_{34}^2 &= -2p_4p_1-p_2p_3, \ &
C_{23}^2 &= (4p_1{}^2+p_2{}^2) L, \ &
C_{13}^2 &=-3Lp_2p_1, \ &
p_4&=0 \ .
\endalignat
$$
The curvature components are
$$
\alignat 1
\Omega_{12} &= K \theta^1 \wedge \theta^2,\
\Omega_{13}  = K \theta^1 \wedge \theta^4 + L(p_1{}^2+4p_2{}^2)\theta^1 \wedge \theta^3-3Lp_2p_1\theta^2  \wedge \theta^3 ,\
\Omega_{14}  = K \theta^1 \wedge \theta^3, \quad \cr
\Omega_{24} &= K \theta^2 \wedge \theta^3,\
\Omega_{23}  = K\theta^2 \wedge \theta^4+L(4p_1{}^2+p_2{}^2) \theta^2\wedge \theta^3 -3Lp_2p_1  \theta^1\wedge \theta^3,\
\Omega_{34} = K \theta^4 \wedge \theta^3. \qquad \quad
\Tag51
\endalignat
$$
The change of basis
$$
\alignat 3
\alpha^1 & = p_1 \theta^2 - p_2 \theta^1 , \quad  &
\alpha^2 & = \theta^3  ,\quad &
\alpha^3 & = -K(\theta^4 -L \theta^3) -p_2 \omega^1 +p_1 \omega^2 ,\\
\alpha^4 & = \sqrt{2}(p_1 \theta^1 +p_2\theta^2) , \quad &
\alpha^5 & = \sqrt{2}(p_1 \omega^1 +p_2\omega^2) , \quad &
\alpha^6 &=  K(\theta^4 +L \theta^3) -p_2 \omega^1 +p_1 \omega^2
\endalignat
$$
for $\lag^*$ leads to the multiplication table {\bf A4} in Theorem 2.2 with isotropy in the dual basis $\lah= span\{ \be_3+\be_6,\be_5\}$.

\smallskip
\noindent
${\text{\bf  Case 4}}$ The subalgebra $F_{14}$ in \cite{pwz:1975a} is $\lah = span \{\bb_1=B_3\} $. By using the basis  $\{\tilde{\bb}_1=B_1, \tilde{\bb}_2=B_2, \tilde{\bb}_3=B_4, \tilde{\bb}_4=B_5, \tilde{\bb}_5=B_6, \bb_1\}$  for $\lo(3,1)$, equations (A.1) and (A.3) give
$$
\alignat 2
\tilde \omega^1 &=  p_1 \theta^1 +p_2 \theta^3, &\quad
\tilde \omega^2 &=  p_3 \theta^1 +p_4 \theta^2 + p_5\theta^3, \quad
\tilde \omega^3  =  p_2 \theta^1 +p_6 \theta^2 + p_7\theta^3 -p_1 \theta^4, \cr
\tilde \omega^4 &= -p_3 \theta^3, &\quad
\tilde \omega^5 &= -p_1 \theta^1. \cr
\endalignat
$$
From condition (A.5), $G/H$ is reductive if and only if $p_4=0$, so we assume $p_4 \neq 0$. The first Bianchi identity gives
$$
C_{14}^1=-p_3{}^2,\, C_{12}^1=-p_2p_3-p_6p_1,\, C_{24}^1=0,\,  C_{23}^1=-3p_7p_3+p_2(p_5-p_6),\,  C_{34}^1=2p_2(p_1+p_4)+2p_5p_3
$$
and
$$
 p_6=t_1 (p_1-p_4), \quad p_5 = t_1(p_4+p_1), \quad
-p_6p_3+p_2p_4-2p_1p_2=0,\quad  p_1 p_3=0
\Tag52
$$
where $t_1\in \reals$.  These last two equations will split into a number of cases. If $p_3 \neq 0$ the Killing form will have rank 4, otherwise the Killing form has rank at most 3, so we split this case into subcases based on $p_3$.

\smallskip
\noindent
${\text{\bf Subcase 4.1}}$ Starting with $p_3 \neq 0$ and $p_1=0$,  we solve \eq{52} and the second Bianchi identity to get
$$
p_2=-t_1 p_3, \quad C_{13}^1=t_2 p_3{}^2,\quad
p_7 =- p_4(2t_1{}^2+t_2)/4
$$
where $t_2 \in \reals$.  The curvature components are
$$
\alignat 1
\Omega_{12} & = K \theta^1\wedge \theta^3 ,
\Omega_{13} = K \theta^1\wedge (\theta^2-t_1\theta^3)
+ p_3{}^2\theta^1 \wedge (t_2\theta^3 -\theta^4)
+\textstyle{ \frac{3}{4}}Lp_3p_4\theta^2\wedge \theta^3,
\Omega_{14} = -p_3{}^2 \theta^1\wedge \theta^3,\qquad
\\
\Omega_{23} &= K\theta^3\wedge( \theta^4 -t_1 \theta^2)+
\textstyle {\frac{1}{2}p_4L(p_4 \theta^2  +\frac{3}{2}p_3 \theta^1)\wedge \theta^3},
\Omega_{24} = 0,
\Omega_{34} = K \theta^2\wedge \theta^3+p_3{}^2\theta^3 \wedge \theta^4
\Tag53
\endalignat
$$
where $K= t_1p_3{}^2$ and $L=t_2-2t_1{}^2$. The change of basis
$$
\eqalign {
\alpha^1 &= -p_3\theta^1 -\textstyle{\frac{1}{2}}p_4\theta^4, \quad
\alpha^2  = p_3\theta^3, \quad
\alpha^4  = -\textstyle{\frac{1}{2}}t_2 p_3 \theta^3-p_3 \theta^4 -\omega^1, \quad
\alpha^5  = -p_4 \theta^2 -t_1p_4 \theta^3, \cr
\alpha^3 &= t_1p_4 \theta^1 +p_3 \theta^4
+t_1(p_4{}^2 -2 p_3{}^2)(2p_3)^{-1} \theta^2
           +(t_1{}^2(p_4{}^2-2p_3{}^2)-t_2 p_3{}^2)(4p_3)^{-1} \theta^3 - \omega^1
}
$$
\noindent
for $\lag^*$ leads to the multiplication table {\bf A1} in Theorem 2.2 with isotropy in the dual basis $\lah= span\{ \be_3 + \be_4 \}$.

\medskip
\noindent
${\text{\bf Subcase 4.2}}$ Starting with $p_3 = 0$ the Bianchi identities give
$$
p_2 = 0,\quad  (p_1-p_4)(C_{13}^1-t_1{}^2p_1p_4) = 0.
$$
The solution to this last equation splits into two further subcases.

\medskip
\noindent
${\text{\bf Subcase 4.2.a:}}$ If  $C_{13}^1 = t_1{} ^2 p_1p_4$ (the Killing form has rank 1), the curvature components are
$$
\alignat 3
\Omega_{12} &= -p_1{}^2 \theta^1 \wedge \theta^2, &\quad
\Omega_{13} &= p_1 L \theta^1\wedge \theta^3 -p_1{}^2 \theta^1 \wedge \theta^4, &\quad
\Omega_{14} &= -p_1{}^2, \\
\Omega_{24} &= -p_1{}^2 \theta^2\wedge \theta^3,  &\quad
\Omega_{23} &= 2L(p_1-p_4)\theta^2\wedge \theta^3 -p_1{}^2 \theta^2\wedge \theta^4, &\quad
\Omega_{34} &= p_1{}^2 \theta^3\wedge \theta^4
\Tag54
\endalignat
$$
where $L = (p_7 +t_1{}^2 p_4)$.  The change of basis
$$
\alpha^1 = \theta^4 -p_7(2 p_4)^{-1} \theta^3, \quad
\alpha^2 = - \theta^2, \quad
\alpha^3 = - \theta^3, \quad
\alpha^5 = \omega^1 -t_1 p_1 \theta^1, \quad
\alpha^6 = -p_4 (\theta^2 + t_1 \theta^3 )
$$
\noindent
for $\lag^*$ leads to the multiplication table {\bf A2} in Theorem 2.2 where $\alpha=p_1/p_4$ and the isotropy in the dual basis is $\lah= span\{ \be_4\}$.

\bigskip
\noindent
${\text{\bf Subcase 4.2.b}}$ We assume $C_{13}^1 - t_1{}^2 p_1p_4 \neq 0 $ so $p_1=p_4 $, and the Killing form has rank 2. Let $L= C_{13}^1+p_4p_7$. The curvature components are
$$
\alignat 3
\Omega_{12} &= -p_4{}^2 \theta^1 \wedge \theta^2 , \quad &
\Omega_{13} &= L\theta^1 \wedge \theta^3 -p_4{}^2 \theta^1 \theta^4 , \quad &
\Omega_{14} &= - p_4{}^2 \theta^1 \wedge \theta^3 , \
\\
\Omega_{23} &= - p_4{}^2 \theta^2 \wedge \theta^4 , \quad &
\Omega_{24} &=  p_4{}^2 \theta^3 \wedge \theta^4 , \quad &
\Omega_{34} &= 0.
\Tag55
\endalignat
$$
\noindent
Write $C_{13}^1 - t_1{}^2 p_4^2 = \epsilon m^2 $ (which is non-zero) where $\epsilon = \pm 1$. The change of basis
$$
\alpha^1 = m \theta^4 -m p_7(2p_4)^{-1} \theta^3 , \quad
\alpha^2 = m \theta^1 , \quad
\alpha^3 = t_1 p_4 \theta^1 -\omega^1 , \quad
\alpha^4 = -p_4\theta^2 -t_1p_4 \theta^3 , \quad
\alpha^5 = m \theta^3
$$
for $\lag^*$ leads to the multiplication table {\bf A3} in Theorem 2.2 with isotropy in the dual basis $\lah= span\{ \be_3 \}$.

\smallskip
\noindent
This concludes the proof of Theorem 2.2. \qed

\bigskip
\noindent
$\underline{\text{\it Proof of Theorem 2.3 - $\lo(2,2)$:}}$  Let $\{ \sigma^i \}_{i=1\ldots 4}$ denote the standard basis for $(\reals^4)^*$, and
$$
\eta^0 =  2 \sigma^1 \sigma^2+
+ 2 \sigma^3  \sigma^4 .
$$
For {\lo}(2,2) use the basis
$$
\alignat 3
2A_1& =e^1_4+e^2_3+e^3_2+e^4_1,\quad & 2A_2 &=e^1_3-e^2_4+e^3_1-e^4_2,\quad &2A_3&=e^1_2-e^2_1+e^3_4-e^4_3 \cr
2B_1&=-e^1_4+e^2_3+e^3_2-e^4_1,\quad &2B_2&=e^1_3+e^2_4+e^3_1+e^4_2,\quad &2B_3&=e^1_2-e^2_1-e^3_4+e^4_3
\endalignat
$$

On pages 2281-2283 of \cite{psw:1977a} the inequivalent subalgebras of $\lo(3,1)$, under inner automorphisms are listed. Of these subalgebras,
labeled $e_{d,n}$, twenty-two  are not reductive in $\lo(2,2)$.

\medskip
\noindent
${\text{\bf Case 1:}}$ We consider the non-reductive subalgebras  of $\lo(2,2)$ which admit a solution to equation (A.5). Therefore, they always lead to a reductive homogeneous space.

\smallskip
\noindent
{\bf Subcase 1.1:} The subalgebra $e_{5,1}$ in \cite{psw:1977a} is $\lah= span\{ \bb_1 = A_2, \bb_2 =A_1-A_3, \bb_3 = B_1, \bb_4 = B_2, \bb_5 = B_3\}$. By using the basis $\{\tilde{\bb}_1 = A_1+A_3, \bb_1, \bb_2, \bb_3, \bb_4, \bb_5 \}$  for {\lo}(2,2), equations (A.1) and (A.3) give $\tilde{\omega}^1=0$. Equation (A.5) has the general solution $\{ r^\alpha_k=0 \}_{ \alpha=1\ldots 5, k=1\ldots 4}$, and $\omega^\alpha \otimes \be_\alpha $ defines a $G$-invariant connection.

\smallskip
\noindent
{\bf Subcase 1.2:} The following 4-dimensional algebras in \cite{psw:1977a} always admit a solution to (A.5). Therefore, they always lead to the construction of a reductive homogeneous space.

$$
\vbox{ \offinterlineskip
\halign{\vrule#\thickspace & $#$\hfil & \thickspace\vrule#\thickspace & \hfil$#$ & \thickspace\vrule#\thickspace & \hfil$#$ & \thickspace\vrule# \cr
\noalign{\hrule}
height2pt&&&&&&\cr
& \lah && \text{Basis\ \ }\{\bb_1,\bb_2,\bb_3,\bb_4\} \hfil && \text{Complement\ \ }\{\tilde{\bb}_1,\tilde{\bb}_2\} \hfil &\cr
height2pt&&&&&&\cr
\noalign{\hrule}
height2pt&&&&&&\cr
& e_{4,2} && A_1-A_3,A_2,B_2,B_1-B_3 && A_1+A_3,B_1+B_3 \hfil &\cr
& e_{4,3} && A_1-A_3,B_1,B_2,B_3 \hfil  && A_1+A_3, A_2\hfil \qquad &\cr
height2pt&&&&&&\cr
\noalign{\hrule}
}
}$$
Equations (A.1) and (A.3) give $\tilde{\omega}^1= \tilde{\omega}^2=0$. Equation (A.5) has the general solution $\{r^\alpha_k=0\}_{\alpha=1\ldots 4, k=1\ldots 4}$, and $\omega^\alpha \otimes \be_\alpha $  defines a $G$-invariant connection.

\smallskip
\noindent
{\bf Subcase 1.3:} The following 3-dimensional algebras in \cite{psw:1977a} always admit a solution to (A.5).  Therefore, they lead to the construction of a reductive homogeneous space.

$$
\vbox{ \offinterlineskip
\halign{\vrule#\thickspace & $#$\hfil & \thickspace\vrule#\thickspace & \hfil$#$ & \thickspace\vrule#\thickspace & \hfil$#$ & \thickspace\vrule# \cr
\noalign{\hrule}
height2pt&&&&&&\cr
& \lah && \text{Basis\ \ }\{\bb_1,\bb_2,\bb_3\} \hfil && \text{Complement \ \ }\{\tilde{\bb}_1,\tilde{\bb}_2\,\tilde{\bb}_3 \} &\cr
height2pt&&&&&&\cr
\noalign{\hrule}
height2pt&&&&&&\cr
&e_{3,2} && B_2,A_2,A_1-A_3 \hfil && A_1+A_3, B_1, B_3 \hfil &\cr
&e_{3,4} && B_3,A_2,A_1-A_3 \hfil && A_1+A_3, B_1, B_2 \hfil &\cr
&e_{3,3};e_{3,5};e_{3,6};e_{3,7}&& A_2+\alpha B_2,A_1-A_3,B_1-B_3 && A_1+A_3, B_1+B_3, B_2 &\cr
height2pt&&&&&&\cr
\noalign{\hrule}
}
}$$
Equations (A.1) and (A.3) give $ \tilde{\omega}^1= \tilde{\omega}^2= \tilde{\omega}^3=0$. Equation (A.5) has the general solution $\{r^\alpha_k=0 \}_{ \alpha=1\ldots 3, k=1\ldots 4}$, and $\omega^\alpha \otimes \be_\alpha $ defines a $G$-invariant connection.

\smallskip
\noindent
{\bf Subcase 1.4:} The following 2-dimensional algebras in \cite{psw:1977a} always admit a solution to (A.5). Therefore, they lead to the construction of a reductive homogeneous space.

$$
\vbox{ \offinterlineskip
\halign{\vrule#\thickspace & $#$\hfil & \thickspace\vrule#\thickspace & \hfil$#$ & \thickspace\vrule#\thickspace & \hfil$#$ & \thickspace\vrule# \cr
\noalign{\hrule}
height2pt&&&&&&\cr
& \lah && \text{Basis\ \ }\{\bb_1,\bb_2\}\hfil && \text{Complement \ \ }\{\tilde{\bb}_1,\tilde{\bb}_2,\tilde{\bb}_3,\tilde{\bb}_4\} &\cr
height2pt&&&&&&\cr
\noalign{\hrule}
height2pt&&&&&&\cr
&e_{2,3}&& B_2, A_1-A_3 \hfil && A_1+A_3, A_2, B_1, B_3 \hfil &\cr
&e_{2,4}&& A_1-A_3, B_3\hfil && A_1+A_3, A_2, B_1, B_2 \hfil &\cr
&e_{2,7}&& A_2, A_1-A_3\hfil && A_1+A_3, B_1, B_2, B_3 \hfil &\cr
&e_{2,10};e_{2,11}&& A_2+c B_2,-A_1+A_3; c>0  \hfil && A_1+A_3, B_1, B_2, B_3 \hfil &\cr
&e_{2,12}&& A_2-c B_3, -A_1+A_3;c \neq 0 \hfil  && A_1+A_3, B_1, B_2, B_3 \hfil &\cr
&e_{2,13}^\epsilon&& B_2+\epsilon(A_3-A_1),B_1-B_3; \epsilon=\pm 1  && A_1+A_3, A_2, B_1+B_3, B_2&\cr
height2pt&&&&&&\cr
\noalign{\hrule}
}
}$$
Equations (A.1) and (A.3) give $\tilde{\omega}^1= \tilde{\omega}^2= \tilde{\omega}^3= \tilde{\omega}^4=0$. Equation (A.5) has the general solution $\{ r^\alpha_k=0 \}_{ \alpha=1\ldots 2, k=1\ldots 4}$, and $\omega^\alpha \otimes \be_\alpha $ defines a $G$-invariant connection.

\smallskip
\noindent
{\bf Subcase 1.4.a:}
The subalgebras $e_{2,8}$ and $e_{2,9}$ in \cite{psw:1977a} are $\lah= span\{ \bb_1 = A_2+B_2, \bb_2=-A_1+A_3+ \epsilon (B_1-B_3) \}$ when $\epsilon=1$ and $\epsilon=-1$ respectively.  By using the basis $\{\tilde{\bb}_1=A_1+A_3, \tilde{\bb}_2=B_1+B_3, \tilde{\bb}_3=A_2-B_3, \tilde{\bb}_3=-A_1+A_3 - \epsilon(B_1-B_3), \bb_1 ,\bb_2 \} $  for {\lo}(2,2), equations (A.1) and (A.3) give
$$
\tilde{\omega}^1=\epsilon(2p_2-p_1)\theta^1, \quad \tilde{\omega}^2=p_1\theta^1, \quad \tilde{\omega}^3=p_2(\theta^4-\epsilon\theta^3), \quad \tilde{\omega}^4=p_2\theta^2.
$$
Equation (A.5) has the general solution
$$
r^1_1=r^1_2=r^2_1=r^2_3=r^2_4=0,\quad r^1_3=\epsilon (p_2-p_1),\quad r^1_4=r^2_2=p_1-p_2,
$$
and $(\omega^1 +(p_1-p_2)(\theta^4-\epsilon\theta^3) )\otimes \be_1 + (\omega^2 + (p_1-p_2) \theta^2 )\otimes \be_2$ defines a $G$-invariant connection.

\medskip
\noindent
We now consider the cases where condition (A.5) is not automatically satisfied.
\smallskip
\noindent
{\bf Case 2:} The subalgebra $e_{2,1}$ in \cite{psw:1977a} is $\lah= span\{A_1-A_3,B_1-B_3\}$.  By using the basis $\{ \tilde{\bb}_1=A_1+A_3, \tilde{\bb}_2=A_2, \tilde{\bb}_3=B_1+B_3, \tilde{\bb}_4=B_2, \bb_1, \bb_2 \}$  for {\lo}(2,2), equations (A.1) and (A.3) give
$$
\tilde{\omega}^1=p_1\theta^1,\quad  \tilde{\omega}^2=p_2\theta^1-2p_1\theta^3,
\tilde{\omega}^3=p_3\theta^1,\quad  \tilde{\omega}^4=p_4\theta^1-2p_3\theta^4.
$$
\noindent
From condition (A.5), $G/H$ is reductive if and only if $p_1=0$ and $p_3=0$. The first Bianchi identity yields
$$
\alignat 5
C_{12}^1 &=\textstyle{\frac{1}{2} (p_2+3p_4)p_3}, &\quad C_{34}^1 &=\textstyle{\frac{3}{2} (p_4-p_2)p_3}, &\quad C_{24}^2 &= -2 p_1 p_3, &\quad C_{23}^2 &=0, &\quad C_{13}^1 &= C_{14}^2,\cr
C_{12}^2 &=\textstyle{\frac{1}{2} (p_4+3p_2)p_1}, &\quad C_{34}^2 &=\textstyle{\frac{3}{2} (p_4-p_2)p_1}, &\quad C_{23}^1 &= -2 p_1 p_3, &\quad C_{24}^1 &= 0.\cr
\endalignat$$
The second Bianchi identity has the general solution
$$p_2 = 0,\quad p_4 = 0, \quad C_{13}^2 = t{p_1}^2,\quad C_{14}^1 = t{p_3}^2,\quad C_{14}^2 = \textstyle{\frac{5}{3} tp_1p_3}$$
where $t\in \reals$.  Let $L=2p_1 p_3$. The curvature components are
$$
\alignat 3
\Omega_{12} &= L\theta^1 \wedge \theta^2, &\quad \Omega_{23} &=-L\theta^1 \wedge \theta^4, &\quad
\Omega_{13} &= \textstyle{ t {p_1}^2\theta^1 \wedge \theta^3 + \frac{5}{6}t L \theta^1 \wedge \theta^4 - L\theta^2 \wedge \theta^4}, \\
\Omega_{34} &= L\theta^3 \wedge \theta^4, &\quad \Omega_{24} &=-L\theta^1 \wedge \theta^3, &\quad
\Omega_{14} &= \textstyle{\frac{5}{6}t L \theta^1 \wedge \theta^3 + t {p_3}^2\theta^1 \wedge \theta^4 - L\theta^2 \wedge \theta^3}.
\Tag56
\endalignat
$$
The Jacobi identities are now satisfied, but depending on the parameters we get non-isomorphic Lie algebras. We now determine the non-isomorphic algebras.

\smallskip
\noindent
{\bf Subcase 2.1:} If $p_1 \neq 0 \text{ and } p_3 \neq 0$ the change of basis
$$
\alignat 3
\alpha^1 &= -p_1 \theta^3 - p_3 \theta^4, &\quad
\alpha^3 &= \textstyle{L(\frac{1}{3}t\theta^1 -\theta^2) - p_1\omega^1 - p_3\omega^2}, &\quad
\alpha^4 &= \sqrt{2}(p_1\theta^3 - p_3\theta^4), \cr
\alpha^2 &= -\theta^1, &\quad
\alpha^6 &= \textstyle{-L(\frac{2}{3}t\theta^1 + \theta^2) + p_1\omega^1 + p_3\omega^2}, &\quad
\alpha^5 &= \sqrt{2}(p_1\omega^5 - p_3\omega^2) \cr
\endalignat
$$
for $\lag^*$ leads to the multiplication table {\bf B2} in Theorem 2.3 with isotropy $\lah= span\{\be_5,\be_3-\be_6\}$.

\smallskip
\noindent
{\bf Subcase 2.2:} If $p_1 = 0 \text{ or } p_3 = 0$ the change of basis
$$
\alignat 3
\alpha^1 &= -p_1 \theta^3 - p_3 \theta^4, &\quad \alpha^3 &= p_1\omega^1 + p_3\omega^2, &\quad
\alpha^5 &= \textstyle{({p_3}^2-{p_1}^2)(\frac{1}{3}t\theta^1 -\theta^2) - p_3\omega^1 + p_1\omega^2,} \cr
\alpha^2 &= \theta^1, &\quad
\alpha^4 &= p_3\theta^3 - p_1\theta^4, &\quad
\alpha^6 &= \textstyle{\frac{1}{2}t(p_3+p_1)^2\theta^2 - ({p_3}^3\omega^1 + {p_1}^3\omega^2)/(p_3+p_1)^2,} \cr
\endalignat
$$
for $\lag^*$ leads to the multiplication table {\bf B3} in Theorem 2.3 with isotropy in the dual basis $\lah= span\{\be_3,\be_5+\be_6\}$ when $p_1=0$ and $\lah= span\{\be_3,\be_5-\be_6\}$ when $p_3=0$. Reversing the sign of $\be_6$ is an automorphism, thus these are equivalent Lie algebra pairs.

\smallskip
\noindent
{\bf Case 3:} The subalgebra $e_{1,10}$ in \cite{psw:1977a} is $\lah= span\{-A_1+A_3\}$.  By using the basis $\{ \tilde{\bb}_1=A_1+A_3, \tilde{\bb}_2=A_2, \tilde{\bb}_3=B_1, \tilde{\bb}_4=B_2, \tilde{\bb}_5=B_3, \bb_1=-A_1+A_3\}$  for {\lo}(2,2), equations (A.1) and (A.3) give
$$
\alignat 2
\tilde{\omega}^1 & = p_1\theta^1 + p_2\theta^4, \quad
\tilde{\omega}^3 = (p_5+p_9)\theta^1 + (p_6-p_3)\theta^4, &\quad
\tilde{\omega}^4 &= (p_3+p_7)\theta^1 + p_8\theta^4,\cr
\tilde{\omega}^2 &= \textstyle{\frac{1}{3}(p_3-p_7)\theta^1 + 2p_2\theta^2 - 2p_1\theta^3 + \frac{1}{3}(p_4-p_9)\theta^4}, &\quad
\tilde{\omega}^5 &= (p_9-p_5)\theta^1 + p_{10}\theta^4.\cr
\endalignat
$$
From condition (A.5), $G/H$ is reductive if and only if $p_1=0$ and $p_2=0$. The first Bianchi identity yields
$$
\alignat 7
p_3&=sp_1, &\quad p_4&=tp_1, &\quad p_5&=-rp_1, &\quad p_8&=tp_1+sp_2, &\quad C_{12}^1&=5p_1p_2J, &\quad C_{13}^1&=-5{p_1}^2J, &\quad C_{23}^1&=-K,\cr
p_9&=sp_2, &\quad p_6&=tp_2, &\quad p_7&=rp_2, &\quad p_{10}&=tp_2+sp_1, &\quad C_{34}^1&=5p_1p_2J, &\quad C_{24}^1&=-5{p_2}^2J \cr
\endalignat
$$
\noindent where $r,s,t \in \reals$ and
$$
J=\textstyle{\frac{1}{3}(rt-s^2),\quad K=\frac{1}{3}({p_1}^2t-2p_1p_2s+{p_2}^2r), \quad L
=\frac{1}{3}({p_1}^2t+4p_1p_2s+{p_2}^2r)\ .}
$$
\noindent The second Bianchi identity has the general solution
$$
C_{14}^1 = 4JK \ .
$$
\noindent The curvature components are
$$
\alignat 1
\Omega_{12} &= p_1(sp_1+rp_2)\theta^1\!\wedge\theta^3\! -\! (K+L)\theta^1\!\wedge\theta^2
\!-\! 5p_1p_2J\theta^1\!\wedge\theta^4 + p_2(tp_1+sp_2)\theta^2\!\wedge\theta^4\! - \! L\theta^3\!\wedge\theta^4, \cr
\Omega_{34} &=  p_1(sp_1+rp_2)\theta^1\!\wedge\theta^3 \! -\! (K+L)\theta^3\!\wedge\theta^4
\!-\! 5p_1p_2J\theta^1\!\wedge\theta^4 + p_2(tp_1+sp_2)\theta^2\!\wedge\theta^4\!  - \! L\theta^1\!\wedge\theta^2, \cr
\Omega_{13} &= p_1(sp_1+rp_2)\theta^1\!\wedge\theta^2\! -\!  2r{p_1}^2\theta^1\!\wedge\theta^3\! +\!  5{p_1}^2J \theta^1\!\wedge\theta^4\! -\! 2sp_1p_2 \theta^2\!\wedge\theta^4\! +\!  p_1(sp_1+rp_2)\theta^3\!\wedge\theta^4, \qquad \cr
\Omega_{24} &= p_2(tp_1+sp_2)\theta^1\!\wedge\theta^2\! -\!  2 t{p_2}^2\theta^2\!\wedge\theta^4 \!  +\!  5{p_2}^2J \theta^1\!\wedge\theta^4\! -\! 2sp_1 p_2 \theta^1\!\wedge\theta^3\! +\!  p_2(tp_1+sp_2)\theta^3\!\wedge\theta^4, \cr
\Omega_{14} &= J(5p_1p_2 \theta^2 - 5{p_1}^2\theta^3 + 4K\theta^4)\!\wedge\theta^1 \! +\!  K\theta^2\!\wedge\theta^3 \! +\!  5J({p_2}^2 \theta^2-p_1p_2\theta^3) \!\wedge\theta^4, \cr
\Omega_{23} &= K\theta^1\!\wedge\theta^4.
\Tag57
\endalignat
$$
\noindent The change of basis
$$
\alignat 1
\alpha^1 &=\textstyle{\frac{1}{3}(2sp_1-rp_2)\theta^1 + p_2\theta^2 - p_1\theta^3 + \frac{1}{3}tp_1\theta^4, \quad \alpha^2 =-p_1\theta^1 - p_2\theta^4,} \cr
\alpha^5 &=\textstyle{\frac{1}{3}(2sp_2+rp_1)\theta^1 - p_1\theta^2 - p_2\theta^3 + \frac{1}{3}tp_2\theta^4, \quad \alpha^4= p_2\theta^1 - p_1\theta^4,} \cr
\alpha^3 &=\textstyle{\frac{1}{9}r(3tp_1-2sp_2)\theta^1 + \frac{1}{3}(tp_1+sp_2)\theta^2 +  \frac{1}{3}(rp_2-3sp_1)\theta^3 + \frac{1}{9}(2stp_1+3trp_2-4s^2p_2)\theta^4 + \omega^1,}
\endalignat
$$
for $\lag^*$ leads to the multiplication table {\bf B1} in Theorem 2.3, with isotropy $\lah= span\{ \be_3\}$.

\smallskip
\noindent
{\bf Case 4:} The subalgebras $e_{1,3}$ and $e_{1,4}$ in \cite{psw:1977a} satisfy Lemma 3.3. Define $A\in GL(4,\reals)$ by $A(f_1)=f_3, A(f_2)=-f_4, A(f_3)=f_1, A(f_4)=-f_2$, where $f_k$ is the standard basis for $\reals^4$, then $A$ is an automorphism of $\lo(2,2)$ that
maps $e_{1,3}$, which has basis $\{ \bb_1= -A_1+A_3-B_1+B_3 \}$, to the
subalgebra $e_{1,4}$, which has basis $\{ \bb_1= -A_1+A_3+B_1-B_3 \}$.
Therefore, we consider only $e_{1,4}$.  By using the basis $\{ \tilde{\bb}_1=A_1-A_3+B_1-B_3, \tilde{\bb}_2=B_2, \tilde{\bb}_3=A_1-A_3+B_1+B_3, \tilde{\bb}_4=-A_1-A_3+B_1-B_3, \tilde{\bb}_5=A_2, \bb_1 \}$  for {\lo}(2,2), equations (A.1) and (A.3) give
$$\alignat 2
\tilde{\omega}^1 &= p_1\theta^1 + p_2\theta^2 + (p_3+p_4)\theta^3 + (p_4-p_3) \theta^4, &\quad
\tilde{\omega}^3 &= (p_7-p_2)\theta^1, \cr
\tilde{\omega}^2 &= 2(p_5-p_3)\theta^1 + (p_6+p_7-p_2)\theta^3 + (p_6-p_7+p_2)\theta^4, &\quad
\tilde{\omega}^4 &= (p_7+p_2)\theta^1, \cr
\tilde{\omega}^5 &= 2(p_5+p_3)\theta^1 + (p_6+p_7+p_2)\theta^3 + (p_6-p_7-p_2)\theta^4. \cr
\endalignat$$

\noindent
From condition (A.5), $G/H$ is reductive if and only if $p_6=0$. The first Bianchi identity yields
$$
\alignat 3
p_4 &=t (p_6+p_2), &\quad C_{24}^1&= {p_7}^2, &\quad C_{34}^1&=2 p_3 p_7+2 p_2 p_4, \quad
C_{12}^1=4 p_2 p_3 - 4 p_3 p_6 - 4 p_5 p_7, \\
p_5 &=t (p_2-p_6), &\quad C_{23}^1&=-{p_7}^2, &\quad C_{14}^1&=-4p_3p_4 + 2{p_2}^2 - 6p_1p_7 - 6{p_7}^2 + 4p_2p_6 + 4p_5p_3 - C_{13}^1 \\
\endalignat
$$
\noindent where $t \in \reals$, and the conditions
$$p_2 p_7=0, \quad   2p_2p_3+(p_3-tp_7)p_6=0.$$

\noindent Let $K=(2p_6t^2 + p_1 + p_7)/p_6$, $L=2(p_6t^2 + p_7 + p_1)/p_6$ and $\beta=C_{13}^1-{p_2}^2 - 2p_2p_6 - 2p_2p_6t^2$. The second Bianchi identity yields $p_3=t p_7$ so the final the remaining conditions are
$$p_2p_7=0,\quad \beta(p_6+p_2)+p_7p_6(L(p_6+2p_7) + Kp_6)=0.$$

\noindent This case splits into three subcases based on $p_7$ and $\beta$.

\smallskip
\noindent
{\bf Subcase 4.1:} If $p_7 \neq 0$ then $p_2=0$. The remaining condition implies $C_{13}^1=-p_7(L(p_6+2p_7) + Kp_6)$. The curvature components are
$$
\alignat 1
\Omega_{23} &=  {p_7}^2\theta^1\wedge ( \theta^4 - \theta^3 ),\quad
\Omega_{12} = -2{p_7}^2\theta^1\wedge ( \theta^2 + t\theta^3 + t\theta^4),\cr
\Omega_{24} &=  {p_7}^2\theta^1\wedge ( \theta^3 - \theta^4 ), \quad
\Omega_{34} = 2t{p_7}^2\theta^1\wedge ( \theta^3 - \theta^4 ), \cr
\Omega_{13} &=-2t{p_7}^2\theta^1\wedge\theta^2 - 2(p_6+p_7)(p_6+2p_7)K\theta^1\wedge\theta^3 - 2({p_6}^2K-{p_7}^2L)\theta^1\wedge\theta^4 \cr
&\quad  - {p_7}^2 \theta^2\wedge ( \theta^3 -\theta^4) + 2t{p_7}^2 \theta^3\wedge\theta^4, \cr
\Omega_{14} &=-2t{p_7}^2\theta^1\wedge\theta^2 - 2(p_6-p_7)(p_6-2p_7)K\theta^1\wedge\theta^4 - 2({p_6}^2K-{p_7}^2L)\theta^1\wedge\theta^3 \cr
&\quad + {p_7}^2 \theta^2\wedge ( \theta^3 -\theta^4) - 2t{p_7}^2 \theta^3\wedge\theta^4,
\Tag58
\endalignat
$$
\noindent The change of basis
$$
\alignat 2
\alpha^1 &= -2 p_6 t\theta^1 + (p_6+2p_7) \theta^3 + (p_6-2p_7) \theta^4, &\quad \alpha^2 &= -4{p_7}^2/p_6\theta^1,\cr
\alpha^3 &= \thickspace \thickspace(p_1+p_7) \theta^1 + p_6 \theta^2 + p_6t \theta^3 + p_6t \theta^4 - p_6/p_7 \omega^1, &\quad \alpha^5 &= -2 p_6 t\theta^1 + p_6 \theta^3 + p_6 \theta^4,\cr
\alpha^4 &= 2(p_1+p_7) \theta^1 - p_6 \theta^2 + p_6t \theta^3 + p_6t \theta^4 - p_6/p_7 \omega^1 \cr
\endalignat$$
for $\lag^*$ leads to the multiplication table {\bf A1} in Theorem 2.3 with isotropy in the dual basis $\lah=span\{\be_3+\be_4\}$.

\smallskip
\noindent
{\bf Subcase 4.2:} If $p_7 = 0$ the remaining the condition splits into two subcases based on $\beta$.

\noindent
{\bf Subcase 4.2.a:} When $\beta=0$, the remaining condition implies $ C_{13}^1=p_2{}^2 + 2p_2 p_6 + 2p_2 p_6 t^2 $. The curvature components are
$$
\alignat 1
\Omega_{12} &= 2{p_2}^2\theta^1 \wedge \theta^2,\quad
\Omega_{23} =-2{p_2}^2\theta^1 \wedge \theta^4,\quad
\Omega_{24} =-2{p_2}^2\theta^1 \wedge \theta^3,\quad
\Omega_{34} = 2{p_2}^2\theta^3 \wedge \theta^4,\\
\Omega_{13} &=-K\theta^1\wedge ( (p_2+2p_6)\theta^3 + (3p_2+2p_6)\theta^4) - 2{p_2}^2\theta^2\wedge\theta^4, \\
\Omega_{14} &=-K\theta^1\wedge ( (p_2+2p_6)\theta^4 + (3p_2+2p_6)\theta^3) - 2{p_2}^2\theta^2\wedge\theta^3.
\Tag59
\endalignat$$

\noindent The change of basis
$$\alignat 1
\alpha^1 &= -p_1 \theta^1 + 2p_6 \theta^2, \quad
\alpha^2 = 2p_6 (\theta^3-\theta^4), \quad
\alpha^3 = 4p_6 \theta^1, \quad
\alpha^4 = p_2 ( \theta^1 + t\theta^3-t\theta^4) +\omega^1, \cr
\alpha^5 &= p_6 ( -2 t\theta^1 + \theta^3 + \theta^4) \cr
\endalignat
$$
for $\lag^*$ leads to the multiplication table {\bf A2} in Theorem 2.3 where $\alpha=-p_2/p_6$ with isotropy in the dual basis $\lah= span\{\be_4\}$.

\smallskip
\noindent
{\bf Subcase 4.2.b:} If $\beta \neq 0$ then $p_2=-p_6$. The curvature components are
$$\alignat 3
\Omega_{12} &= 2{p_6}^2\theta^1\!\wedge\theta^2, &\quad
\Omega_{23} &=-2{p_6}^2\theta^1\!\wedge\theta^4, &\quad
\Omega_{13} &= (\beta -{p_6}^2K)\theta^1\!\wedge (\theta^3 - \theta^4) - 2{p_6}^2 \theta^2\!\wedge\theta^4,\\
\Omega_{34} &= 2{p_6}^2\theta^3\!\wedge\theta^4, &\quad
\Omega_{24} &=-2{p_6}^2\theta^1\!\wedge\theta^3, &\quad
\Omega_{14} &= (\beta -{p_6}^2K)\theta^1\!\wedge ( \theta^4 - \theta^3) - 2{p_6}^2 \theta^2\!\wedge\theta^3.\quad\quad
\Tag60
\endalignat
$$
\noindent The change of basis
$$\halign{\quad$#$\hfil &\quad$#$\hfil &\quad$#$\hfil \cr
\alpha^1 = \sqrt{|\beta/2|} ( \frac{1}{2}p_1/p_6 \theta^1 -\theta^2), &
\alpha^3 = p_6 ( \theta^1 + t\theta^3 - t\theta^4) - \omega^1, &
\alpha^5 =\sqrt{|2\beta|} \theta^1, \cr
\alpha^2 =\sqrt{|\beta/2|} (\theta^3 - \theta^4), &
\alpha^4 = p_6 ( -2t\theta^1 + \theta^3 + \theta^4), \cr
}
$$
for $\lag^*$ leads to the multiplication table {\bf A3} in Theorem 2.3 where $\epsilon = \frac{\beta}{|\beta|}$ and the isotropy in the dual basis is $\lah= span\{\be_3\}$.

This concludes the proof of Theorem 2.3.\qed

We now list the algebra pairs in Theorems 2.2 and 2.3 where the metrics can be Einstein without being of constant curvature.

\proclaim{Lemma 5.1} Let $(G/H,\eta)$ be a homogeneous non-reductive pseudo-Riemannian Einstein manifold of dimension four which is not of constant curvature
and where $H$ is connected.
\roster
\item"{\rm i]}" The space is Ricci-flat if and only if the Lie algebra pair $(\lag, \lah)$ is isomorphic to {\bf B3}.

\item"{\rm ii]}" The space is Einstein and not Ricci flat if and only if the Lie algebra is isomorphic to the pair {\bf A2} with $\alpha = 2/3$.
\endroster
\endproclaim
\demo{Proof} Starting with the curvature forms in equations
\equationlabel{5}{51}{}, \equationlabel{5}{53}{}, \equationlabel{5}{54}{}, \equationlabel{5}{55}{},
\equationlabel{5}{56}{}, \equationlabel{5}{57}{}, \equationlabel{5}{58}{}, \equationlabel{5}{59}{}, \equationlabel{5}{60}{}, the coefficients of the Ricci tensor in an orthonormal frame $u_\sigma$ at $\sigma = [H]$ are easily computed. For example, from the curvature in \equationlabel{5}{54}{} we get
$$
Ricci = (3 \alpha^2 p_4{}^2) \ \eta^0  -p_4(t_1{}^2 p_4+p_7)(3\alpha -2) \theta^3 \otimes \theta^3\,
$$
where $p_4 \neq 0$. If $\alpha = \frac{2}{3}$ and $t_1{}^2p_4 +p_4 \neq 0$ the
space is Einstein and not of constant curvature. Similar computations with \equationlabel{5}{59}{} and \equationlabel{5}{56}{} in {\bf Subcase 2.2} prove the lemma. \qed
\enddemo

Case i] in Lemma 5.1 corresponds to Proposition 2.5.2 on pg. 153 of \cite{komrakov:2001a}, and
case ii] in Lemma 5.1 corresponds to Proposition 1.4.2 on pg. 142 of \cite{komrakov:2001a}.

\initializesection6

\head 6. Global Results and Existence \endhead

To prove Theorem 2.4, we start by characterizing the four-dimensional simply connected non-reductive pseudo-Riemannian homogeneous spaces. These turn out to be fairly simple.

\proclaim{Theorem 6.1} Let $(\lag,\lah)$ be a Lie algebra pair from
Theorem 2.2 or Theorem 2.3 and suppose $G$ is the simply connected Lie group with Lie algebra $\lag$. Then there exists a closed connected Lie subgroup $H \subset G$ with Lie algebra $\lah$ such that $G/H$ is diffeomorphic to $\reals^4$.
\endproclaim
\demo{Proof} The proof is done on a case by case basis. We first
consider the pairs in {\bf A2}-{\bf A3} of Theorem 2.2 (or Theorem
2.3) where $\lag$ is solvable. Let $G$ be the simply connected
solvable Lie group having Lie algebra $\lag$, and $H$ be the
connected Lie subgroup having Lie subalgebra $\lah$. Since $H$ is
closed \cite{chevalley:1969a}, $G/H$ is a manifold. Since $H$ is
connected, $G/H$ is simply connected and $G/H$ is diffeomorphic to
$\reals^4$ \cite{mostow:1954a}.

For cases {\bf A1} and {\bf A4} of Theorem 2.2 and {\bf B1-B3} of Theorem 2.3,  we construct connected Lie groups $G^0$ and connected closed subgroups $H^0 \subset G^0$ such that the covering space of $G^0/H^0$ is $\reals^4$.
It follows
that $\reals^4 = \tilde G/H$ (see Theorem 2.1 p. 125 \cite{onishichik:1998a}) where $\tilde G$ is the simply connected cover of $G^0$ and $H$ is a closed connected Lie subgroup having Lie subalgebra $\lah$. We start with {\bf A4} and {\bf B2}.

Let $\ba, \bb \in \reals^2$, and $\ba \times \bb = a_1  b_2 -b_1 a_2$.  The multiplication map for the six-dimensional Lie group $SL(2,\reals) \sd N_3 $ is
$$
\left(A,\ba ,\alpha \right) *
\left(B,\bb,\beta \right)
= ( A A' ,   A \bb + \ba , \alpha+\beta -  (A \bb) \times \ba)
$$
where $ A,B \in SL(2,\reals), \ba,\bb \in \reals^2$ and $\alpha,\beta \in \reals$. Let $H^0_l$ and $H^0_n$ be the closed subgroups
$$
H^0_l= \{ \left( \matrix 1 & 0 \\ t & 1 \endmatrix \right) , \left( \matrix 0 \\ s \endmatrix \right),  2 t )  \quad t,s \in \reals \} \ , \quad
H^0_n= \{ \left( \matrix 1 & 0 \\ t & 1 \endmatrix \right) , \left( \matrix 0 \\ s \endmatrix \right),  -2 t )  \quad t,s \in \reals \} \ .
$$
The Lie algebra pair in {\bf A4} is isomorphic to $(\lag,\lah_l)$, and the pair in {\bf B2} is isomorphic to $(\lag,\lah_n)$. The quotient spaces $G^0/H^0_l$ and $G^0/H^0_n$ are diffeomorphic to $(\reals^2\backslash \{(0,0)\}) \times \reals^2 $, so the covering space in these cases is $\reals ^4$.

For the Lie algebra pair {\bf B1} in Theorem 2.3, let $G^0$ be the group $ SL(2,\reals)\sd \reals^2 $, and let
$$
H^0 =
\{ \left( \matrix 1 & 0 \\ t & 1 \endmatrix \right) , \left( \matrix 0 \\ 0 \endmatrix \right)  \quad t \in \reals \}\ .
$$
The pair in {\bf B1} of Theorem 2.3 is isomorphic to this $(\lag,\lah)$.
The quotient space $G^0/H^0$ is diffeomorphic to $(\reals^2\backslash \{(0,0)\})  \times \reals^2$, so its simply connected cover is $\reals ^4$.
The Lie algebra pair in {\bf A1} in Theorems 2.2 and 2.3 is similar to this one.

For \A5, consider the monomorphism $\phi:\lag \to \lo(2,3)$ given by
$$
\alignat 2
\phi(\be_1) & = e^1_4+e^4_1-e^2_3-e^3_2,& \quad
\phi(2 \be_2) &= e^1_2+e^1_3+e^2_4 +e^3_1+e^4_2+e^4_3-e^2_1-e^3_4, \\
\phi(\sqrt{2}\be_5) &= e^1_5 +e^5_1+e^5_4-e^4_5,& \
\phi(2 \be_3) &= e^1_3+e^2_1 +e^2_4+e^3_1+e^3_4+e^4_2-e^1_2-e^4_3, \\
\phi(\sqrt{2}\be_6) &=  e^3_5-e^2_5-e^5_2-e^5_3, & \
\phi(2\be_4) &= e^2_1+e^2_4+e^4_2+e^4_3-e^1_2-e^1_3-e^3_1-e^3_4, \\
\phi(\be_7) &= e^1_4+e^2_3+e^3_2+e^4_1,& &
\endalignat
$$
where $\eta^0 =diag(-1,-1,1,1,1)$ is the symmetric $5\times 5$ matrix defining $\lo(2,3)$. Let $G$ be the simply connected seven-dimensional Lie group having Lie algebra $\lag$, and let $\Phi:G \to O(2,3)$ be the induced homomorphism from $\phi$. We now show that $G$ acts transitively on the manifold
$$
M= \{ (x_1,x_2,x_3,x_4,x_5 )\in \reals^5 \ | \ -x_1{}^2-x_2{}^2+x_3{}^2+x_4{}^2+x_5{}^2=-r^2 \}
$$
which is diffeomorphic to $S^1 \times \reals^3$ \cite{wolf:1984a}. Let $(x_1,x_2,x_3,x_4,x_5) \in M$. Application of the group element $e^{\phi(t\be_3-t\be_2)}$, where $t=\pi/2$ if $x_2=0$ otherwise $\tan t = x_1/x_2$, maps this point to $(0,\tilde x_2,\tilde x_3,\tilde x_4,\tilde x_5)$, where $\tilde x_2\neq 0$.
Similar use of the one parameter subgroups of $G$ map this point to $(0,r,0,0,0)$. Hence, $G$ acts transitively on $ M$. The Lie algebra of the Lie subgroup of $G$ which stabilizes $(0,r,0,0,0)$ is $\lah=\{ \be_1+\be_7, \be_3-\be_4,\be_5\}$. Therefore, the covering space for $G/H$ is $\reals^4$.

The details for {\bf B3} can be found in the proof of Theorem 2.5 in
section 7.  \qed
\enddemo

Theorem 6.1  can now be used to prove Theorem 2.4 i].

\demo{Proof} Proof of Theorem 2.4 i]: Since the homogeneous space $G/H$ in the theorem is simply connected, we may assume $G$ is simply connected and $H$ is connected. By Theorem 2.2 or Theorem 2.3,
the Lie algebra pair ($\lag$,$\lah$) is isomorphic to one of {\bf A1-}\A5 or {\bf B1-B3}. The Lie algebra isomorphism lifts to a Lie group isomorphism to one of the simply connected groups used in the proof of Theorem 6.1. Therefore, $M$ is diffeomorphic to $\reals^4$.
\enddemo

Up to this point, we have shown that a simply connected non-reductive homogeneous pseudo-Riemannian manifold is diffeomorphic to $ \reals^4$ and the Lie algebra of its isometry group must be isomorphic to one in Theorems 2.2 or 2.3. We now show that \A5 cannot occur.

\proclaim{Lemma 6.1} Let $(G/H,\eta)$ be a simply connected four-dimensional homogeneous Lorentz manifold with Lie algebra pair \A5.  Then $G$ is a proper subgroup of the isometry group $\tilde O(2,3)$.
\endproclaim

\demo{Proof} The computations in Theorem 6.1 show that there exists a transitive action of $G$ on $S^1\times \reals^3$ with isotropy $K$ which has the same Lie algebra as $H$.
We showed in the proof of Theorem 2.2 in section 5 {\bf Case 2}, that the Lorentz metric, which is unique up to scaling, was of constant (negative) curvature,
so the standard action of $O(2,3)$ on $S^1\times \reals^3$
 is by isometries for this metric. Therefore, an invariant Lorentz metric on $G/H$ (the covering
space) will admit $\tilde O(2,3)$ acting by isometries, and the Lie algebra of the isometry group will not be $\lag$.  \qed
\enddemo

Lemma 6.1 allows us to prove Theorem 2.4 ii] by eliminating \A5.

\demo{Proof} Proof of Theorem 2.4 ii]: Starting with $G/H$ simply
connected we may assume that $H$ is connected and so Theorem 2.2, or
2.3 imply that the Lie algebra pair is isomorphic to one in the
lists in these two theorems. However if $G$ is the isometry  group
then by Lemma 6.1, \A5 cannot be the Lie algebra pair for the
isometry group of simply connected four-dimensional homogeneous
Lorentz manifold. \qed
\enddemo

In order to finish the proof of Theorem 2.4 (the converse part) we need to show that we can build metrics on $ \reals^4$ having the isometry algebras in Theorems 2.2 (except \A5) and 2.3. In order to do this we first give two lemmas.

\proclaim{Lemma 6.2} Let $G/H$ be a homogeneous space with pair $(\lag,\lah)$, and $H$ connected. If the pair $(\lag,\lah)$ satisfies Lemma 3.2, then
$$\
\eta(X,Y) = \eta^0(\theta(X),\theta(Y)) \qquad X,Y \in T_p G
$$
is basic for the projection $\bq :G \to G/H$ and defines a pseudo-Riemannian metric on $G/H$ with curvature tensor $\Omega$.
\endproclaim
\demo{Proof} The form $\theta$ is $\reals^n$ valued, so $\eta$
defines a symmetric bilinear form on $TG$. From $1]$ in Lemma 3.1, the
form $\eta $ is semi-basic for the projection $\bq: G \to G/H$ and
has the same signature as $\eta^0$. The Lie derivative of $\eta$ with respect to $X\in \lah$,
$$
{\Cal L}_X \eta^0( \theta , \theta) = \eta^0( {\Cal L}_X \theta,\theta) + \eta
(\theta, {\Cal L}_X \theta) = \eta^0 (\rho_*(X) \theta,\theta) +\eta^0(\theta,
\rho_*(X) \theta) =0 ,
$$
implies that $\eta$ is $H$-basic, because $H$ is connected.
We can check that $\Omega$ is the curvature of $\eta$,
by choosing a local cross section of $\bq:G \to G/H$ and pulling
back the structure equations by the cross section, or by
reversing the arguments in section 3 which we now do.
Let $u_\sigma \in O(M)$, and $\Psi$ be constructed
as in \equationlabel{3}{33}{}. The pullback $\Psi^* {\bar \theta}$ of the canonical form $\bar \theta$ on the frame bundle, are $G$-invariant and provide a basis for the $\bq :G \to G/H$ semi-basic forms. Therefore,
$$
\theta = A \Psi^*{\bar \theta} \quad {\text {\rm where}}\ \ A \in GL(n,\reals).
\Tag61
$$
Now, let $X,Y\in \reals^n$ and choose $\tilde X, \tilde Y \in T_eG$ such that
$\bq_* \tilde X  = u_\sigma X$, and
$\bq_* \tilde Y  = u_\sigma Y$. By definition of $u_\sigma$ and $\eta$, we
have
$$
\eta^0(X,Y) = \eta( u_\sigma X, u_\sigma Y) = \eta^0(\theta (\tilde X),\theta(
\tilde Y)) \, .
$$
Condition \eq{61} gives
$$
\eta^0(X,Y) = \eta^0(A {\bar \theta}(\Psi_* X),A {\bar \theta}(\Psi_*Y)) =
\eta^0(A u_\sigma ^{-1}\pi_* \Psi_* X, A  u_\sigma ^{-1}\pi_* \Psi_* X )\, .
$$
The commutative diagram \equationlabel{3}{34}{} gives
$$
\eta^0(X,Y) = \eta^0(A \bq_* (\tilde X) , A \bq_*(\tilde Y) ) = \eta^0(AX, AY),
$$
so $A \in O(p,\tilde p)$. Finally, using the frame
$v_\sigma = u_\sigma A$ to redefine $\Psi$, we get $\Psi^* \theta = \theta , \Psi^* \omega = \omega $, and $\Psi^* \Omega = \Omega$. \qed
\enddemo

This lemma says that for any case we consider in section 5, and no matter what value we choose for the parameters in the curvature form $\Omega$, we
can construct a homogeneous pseudo-Riemannian manifold having the chosen value
of the curvature form. In the next lemma, we give a sufficienct condition
on the curvature for a given Lie algebra to be the Lie algebra of
the isometry group.

\proclaim{Lemma 6.3} Let $(G/H,\eta)$ be an $n$-dimensional pseudo-Riemannian homogeneous space with curvature form $\Omega$ and let $R_{ijkl}$ and $R_{ijkl;m}$ be the corresponding components of the Riemann curvature tensor and its covariant derivative in the orthonormal frame $u_\sigma$ at $\sigma =[H]$. Let
$$
\eqalign {
S  = \{   E \in \lo(p,\tilde p) \ | & R_{sjkl} E^s_{i} +  R_{iskl} E^s_{j} + R_{ijsl} E^s_{k} + R_{ijks} E^s_{l} = 0, \cr
 & R_{sjkl;m} E^s_{i} +  R_{iskl;m} E^s_{j} + R_{ijsl;m} E^s_{k}
 + R_{ijks;m} E^s_{l} +  R_{ijkl;s} E^s_{m} = 0
\} .
}
\Tag62
$$
If $\dim S= \dim\lah $, then the Lie algebra of the isometry group is $\lag$ (the Lie algebra of $G$).
\endproclaim

\demo{Proof} Use the notation in section 3. If $E = \rho_*(\be)$ where $\be \in \lah$, then $E \in S$, so $ \dim S \geq \dim \lah $. Suppose that $\tilde G$ is the isometry group of $(G/H,\eta)$. To prove the lemma, it  is sufficient to show $\dim \tilde G = \dim G$. Let
$\tilde H \subset \tilde G$ be the isotropy subgroup at $\sigma = [H] \in G/H$
with linear isotropy representation $\tilde{\rho}$. We have $G \subset \tilde G$, $ \dim G = n + \dim H$, and $\dim \tilde G = n + \dim \tilde H$. By the argument just given, $\tilde{\rho}_* (\tilde {\lah})$ satisfies \eq{62}. Therefore, if the hypothesis of the theorem hold then $\dim \tilde {\lah} = \dim \lah $ and $ \dim G = \dim \tilde G$. \qed
\enddemo

\noindent
The set $S$ is the Lie algebra of the stabilizer of both the curvature tensor and its covariant derivative at a point. This lemma states that if this subalgebra has the same dimension as $\lah$, then the isometry algebra can not have dimension greater than $n+\dim \lah$. Therefore, it must be  the given algebra. Lemma 6.3 of course generalizes using the higher order covariant derivatives of the curvature tensor.

We can now prove the converse condition in Theorem 2.4 by using Lemmas 6.2,  Lemma 6.3, and the computations from section 5. That is, for each Lie algebra pair in the theorem we find values for the coefficients of the curvature form $\Omega$ such that Lemma 6.3 is satisfied.

\demo{Proof} Proof of converse for Theorem 2.4: We start with the Lorentz signature.

\noindent {\bf Case 3:} If $C_{13}^1(p_1{}^2+4p_2{}^2)^{-1} =L \neq
0$, then Lemma 6.3 is satisfied, and the isometry algebra is {\bf
A4}.

\noindent
{\bf Subcase 4.1:} If $ t_2-2t_1{}^2 =L  \neq 0 $ then Lemma 6.3 is satisfied, and the algebra is {\bf A1}.

\noindent
{\bf Subcase 4.2.a:} If $p_7+t^2p_4 \neq 0$ then Lemma 6.3 is satisfied, and the algebra is {\bf A2}.

\noindent
{\bf Subcase 4.2.b:}  If $L=C_{13}^1+p_7p_4 \neq 0 $ then Lemma 6.3 is satisfied, and the algebra is {\bf A3}.

\medskip
\noindent
Now we consider the signature $(2,2)$ cases.
\smallskip

\noindent
{\bf Subcase 2.1:}  If $t \neq 0 $ then Lemma 6.3 is satisfied, and the isometry algebra is {\bf B2}.

\noindent {\bf Subcase 2.2:}  See the first part of the of Theorem
2.5 in section 7.

\noindent
{\bf Case 3:}  If $rt-s^2\neq 0 $ then Lemma 6.3 is satisfied, and the isometry algebra is {\bf B1}.

\noindent {\bf Subcase 4.1:}  If $2p_6t^2+p_1+p_7 =K \neq 0$, Lemma
6.3 is satisfied and the isometry algebra is {\bf A1}.

\noindent {\bf Subcase 4.2:}  If $p_2\neq 0 $ and $ 2p_6 t^2+p_1
\neq 0 $, Lemma 6.3 is satisfied and the isometry algebra is {\bf
A2}.

\noindent {\bf Subcase 4.3:}  If $C^1_{13}-p_1p_6+p_6^2 \neq 0 $,
Lemma 6.3 is satisfied, and the isometry algebra is {\bf A3}. \qed
\enddemo

\initializesection7

\head 7.  The Einstein Examples \endhead

In this section, we prove Theorem 2.5 by  constructing all
homogeneous Einstein and Ricci flat metrics on the simply connected
non-reductive homogeneous spaces of dimension four.

\demo{Proof} (Theorem 2.5) By Lemma 5.1, the Lie algebra pair
$(\lag,\lah)$ of $G/H$ is isomorphic to {\bf B3} if $\eta$ is
Ricci-flat, otherwise it is isomorphic to the Lie algebra pair {\bf
A2} with $\alpha=\frac{2}{3}$. To prove the theorem, it is
sufficient to construct the two simply connected homogeneous spaces
that have Lie algebra pair {\bf B3} or {\bf A3} with $\alpha =
\frac{2}{3}$, and find all the invariant metrics. Theorem 2.4 says
that the manifolds themselves are diffeomorphic to $\reals^4$.

We remind the reader that the Lie algebra of infinitesimal
generators of $G$ acting on $G/H$  is isomorphic to the Lie algebra
of right invariant vector fields, and so we use a basis of left
invariant forms $\sigma$ which have structure constants negative to
the ones in Theorem 2.2 or 2.3 to construct our examples.

\smallskip

We start by proving the second part of the theorem using the Lie
algebra in {\bf A2}.  The (negative of the) structure equations are
easily integrated on $\reals^5$ to give the left-invariant forms
$$
\sigma^5 = dx^5, \ \sigma^4 = e^{-x^5} dx^4,\
\sigma^3=e^{(1-\alpha)x^5} dx^3, \ \sigma^2 = e^{-\alpha x^5}(dx^2
+y^3 dx^4) , \ \sigma^1 = e^{-(\alpha+1)x^5}(dx^1-x^4 dx^2) \ .
$$
Let $X_i, i=1\ldots 5$ be the dual frame of left invariant vector
fields. The $X_4 $ basic symmetric bilinear forms on $G$ are easily
computed to be
$$
\tilde \eta = a_1 (2\sigma^1 \sigma ^3 - (\sigma^2)^2  ) + a_2
(\sigma^3)^2  +2 a_3 \sigma^3  \sigma^5 + a_4 (\sigma^5)^2 .
$$
With the coordinates $y^1=x^1 +x^3(x^4)^2/2,y^2=x^2
+x^3x^4,y^3=x^3,$ and $ y^4=x^5$ on the quotient of $G$ by the
orbits of $X_4 =
e^{x^5}(\partial_{x^4}-x^4x^3\partial_{x^1}-x^3\partial_{x^2})$ we
have $\tilde \eta = \pi^* \eta$ where
$$
\eta = a_1 e^{-2\alpha y^4}(2 dy^1 dy ^3 - (dy^2)^2 ) + a_2
e^{2(1-\alpha)y^4} (dy^3)^2  +2 a_3 e^{(1-\alpha)y^4} dy^3 dy^4 +
a_4 (dy^4)^2  .
$$
These bi-linear forms are non-degenerate if $a_1a_4 \neq 0$. If $a_2
\neq 0$ the only  Killing vectors are
$$
Y_1 = \partial_{y^1},\ Y_2 = \partial_{y^2},\ Y_3 = \partial_{y^3},\
Y_4 = y^2\partial_{y^1}+ y^3\partial_{y^2} , \ Y_5 =
(1+\alpha)y^1\partial_{y^1} +\alpha y^2\partial_{y^2}+(\alpha -1)
y^3\partial_{y^3} +\partial _{y^4} .
$$
These vector fields form a Lie algebra with the multiplication table
in {\bf A2}. At the point $(0,0,0,0)$ the isotropy is $Y_4$, and
$\eta$ is the most general metric invariant under the flow of these
Killing fields. When $\alpha = \frac{2}{3}$, these metrics are
Einstein and they are not of constant curvature when $a_2 \neq 0$.
The signature of theses metrics can be only Lorentz or $(2,2)$, as
expected.

\smallskip
For the Lie algebra pair in {\bf B3}, consider the quotient
 space $G/H$ where $G=SL(2,\reals) \sd \reals^2 \oplus \reals$ and
$$
H= \{ \left( \matrix 1 & 0 \\ t & 1 \endmatrix \right) , \left( \matrix 0 \\
s \endmatrix \right),   s )  \quad t,s \in \reals \}  \ .
$$
The manifold $G/H$ is diffeomorphic to $\reals^2-{(0,0)} \times \reals^2$. In
 terms of coordinates $(x^1,x^2,x^3,x^4)$ on $G/H$ the most general metric having
  signature $(2,2)$  on $G/H$ where $G$ acts by isometries is
$$
\hat \eta =  2 \sigma^1  \sigma^2 + 2 \sigma^3  \sigma^4
$$
where $\sigma^1=dx^1,\sigma^2 =
dx^4+\frac{Lx^2}{2}(x^2dx^1-x^1dx^2), \sigma^3 = dx^2,
\sigma^4=-dx^3 +\frac{Lx^1}{2}(x^1 dx^2-x^2dx^1) $. Let $ \pi
:\reals^4 \to M$ be given by $ x^1= e^{y^1}\cos y^2, x^2=e^{y^2}
\sin y^2, x^3=y^3, x^4=y^4 $, and $\eta = \pi^* \hat \eta$. This
gives the metrics $\eta$ in part one of Theorem 2.5. The covering
group $\tilde G = \tilde {SL}(2,\reals) \sd \reals^2 \oplus \reals $
  acts transitively and by isometries on $\reals^4$ for the metrics $ \eta$.
 There are eight Killing vector fields for $ \eta $:
$$
\eqalign  {
Y_1& =\cos(2y^2)\partial_{y^1}- \sin(2y^2) \partial_{y^2} +y^3\partial_{y^3}-
y^4\partial_{y^4}  , \
Y_2 =\textstyle{ \frac{1}{2}\sin(2y^2)\partial_{y^1}+\cos^2 y^2\partial_{y^2}
 +y^3\partial_{y^4}\,} \cr
Y_3 &=\textstyle{\frac{1}{2}\sin(2y^2)\partial_{y^1}\!-\!\sin^2 (y^2)\partial
_{y^2} \!+\!y^4 \partial _{y^3}}, \
Y_4 =\partial_{y^4}, \
Y_5  =-\partial_{y^3} , \
Y_6 =e^{y^1}\cos(y^2)\partial_{y^3}\!+\!e^{y^1}\sin y^2\partial_{y^4}, \cr
Y_7 &=e^{-y^1}\cos(y^2) \partial_{s^1}\!-\!e^{y^1}\sin(y^2)\partial_{s^2} \!+
\!
Le^{2y^1}(\textstyle{\frac{1}{2}}\sin(2y^2)\partial_{y^3}\!+\!\sin^2(y^2) \partial_{y^4})\, , \  \cr
Y_8 &=e^{-y^1}\sin(y^2) \partial_{s^1}+e^{y^1}\cos y^2\partial_{s^2} -
Le^{2y^1}(\cos^2 (y^2) \partial_{y^3}+\textstyle{\frac{1}{2}}\sin(2y^2) \partial_{y^4})\, .  \
}
$$
The first six of these vector fields are complete and are the
infinitesimal generators for the transitive action of $\tilde G$ on
$\reals^4$. The metric $\eta$ is the most general metric invariant
under the flow of these six Killing fields. In order to have
 a complete Killing vector field of the form $ a X_7 +b X_8$, it is
necessary and sufficient that the differential equations
$$
\frac{dy^1}{dt} =  e^{-y^1}\left(a  \cos y^2+ b \sin y^2 \right) ,
\frac{dy^2}{dt} =  e^{-y^1}\left(b  \cos y^2- a \sin y^2 \right)
$$
admit solutions defined for all $t\in \reals$. By letting $z=y^1-i y^2$, this
equation becomes
$$
\frac{dz}{dt} = (a+i b ) e^{-z}
$$
which has solution $ z = \ln \left( (a+ib) t+(c_0+i c_1)\right)$.
For any initial condition of the form $ z(0) = r (a+ib) $ where $
r\in \reals$, the solution does not exist for all $t  \in \reals$.
Therefore, the Lie algebra of the isometry group is
$\Sl(2,\reals)\sd\reals^2\oplus\reals$, and at the point $(0,0,0,0)$
the isotropy subalgebra is $\{Y_3,Y_5+Y_6\}$.  This proves the first
part of the theorem. For completeness, the non-zero components of
the curvature form are
$$
\Omega^1_3 =-\Omega^3_2 = - 3 L {\hat \sigma}^1 \wedge {\hat \sigma}^3\ ,
$$
where ${\hat\sigma}^i = \pi^* \sigma^i $. By using this formula for
the curvature form, the Ricci tensor is found to vanish. It is worth
noting, that the covariant derivative of the Riemann curvature
vanishes. \qed
\enddemo

\initializesection8

\head 8.  Concluding remarks \endhead

In this paper we have investigated only the most basic questions surrounding the geometry
of non-reductive homogeneous spaces. We have not addressed such natural problems as determining the  holonomy of these spaces \cite{ghanam:2001a}, the homogeneous structure of the geodesics of theses spaces \cite{kowalski:2000a}, or whether they are geodesically complete. These problems will be considered in future research.

\initializesection9

\head Appendix A. Some Computations in a basis \endhead

We write out the conditions in Lemma 3.2 in a basis. Let $ \{
\bb_\alpha\}_{ \alpha = 1\ldots q}$ be a basis for the subalgebra $
\lah \subset \lo(p,\tilde p)$ and complete this to a basis $ \{
\tilde {\bb} _r,\bb_\alpha \}_{ r=1\ldots n(n-1)/2-q,\alpha =
1\ldots q} $ for $\lo(p,\tilde p)$. In this basis, the structure
constants are
$$
[\bb_\alpha,\bb_\beta] = K^\gamma_{\alpha \beta} \bb_ \gamma  \ , \qquad
[\bb_\alpha,\bb_r] = \tilde K^\beta_{\alpha r} \bb_ \beta + \hat K^s_{\alpha r} \bb_ s \ , \qquad
[\bb_r,\bb_s] = \tilde L^\alpha_{r s} \bb_\alpha +  L^t_{r s} \bb_t  \ .
$$
Let $\{\be_\alpha \}_{\alpha=1\dots q}$, form a basis for $\lah
\subset \lag = \lah \oplus \reals^n$, where $\rho_*(\be_\alpha)=
\bb_\alpha$, and complete this to a basis $\{\tilde{\be}_i,
\be_\alpha\}_ {1\leq i \leq n, 1 \leq \alpha \leq q}$ for $\lag$.
Let $ \{\theta^i, \omega^\alpha\}_{i=1\ldots n,\alpha=1\ldots q} $
be the dual basis. We may then write
$$
\omega = \omega^\alpha \bb_\alpha+ \tilde \omega^r \tilde {\bb}_r \
, \quad \Omega = \Omega^\alpha \bb_\alpha+ \tilde \Omega ^r \tilde
{\bb}_r
$$
where $\omega^\alpha,\tilde \omega^r, \Omega^\alpha,$ and $\tilde
\Omega ^r  \in \lag ^*$. By conditions $\hat 2]$ and $\hat 3]$,
these forms satisfy
$$
\omega^\alpha( \be_\beta  ) =\delta ^\alpha_\beta\, ,\quad \tilde
\omega^r(\be_\alpha) = 0 \, , \quad \Omega^\alpha(\be_\alpha) = 0 \,
, \quad \tilde \Omega^r(\be_\alpha) = 0\, .
$$
Consequently, we may write
$$
{\tilde \omega}^r = P^r _i \theta^i , \quad \Omega^\alpha =
\Omega^\alpha_{ij} \theta^i \wedge \theta^j, \quad {\tilde \Omega}^r
= {\tilde \Omega}^r_{ij} \theta^i \wedge \theta^j \tag{A.1}
$$
where $ P^r_i, \Omega^\alpha_{ij},$ and ${\tilde \Omega}^r_{ij} \in
\reals$. By using $\hat 2]$ and $\hat 3] $,  we also have
$$
\iota_{\be_\alpha} d\theta^i = - b_{\beta j}^i \omega^\beta(\be_\alpha) \wedge \theta^j  = - b^i_{\alpha j }\theta^j
$$
where $\bb_\alpha = (b_{\alpha j}^i) $.  From this equation, it follows
$$
d\theta^i = -  b^i_{\alpha j} \omega ^ \alpha \wedge \theta ^j - b_{r [k}^i
 P^r_{j]} \theta^j \wedge \theta ^k \ .
\tag{A.2}
$$
By (A.1), $\omega = \omega^\alpha \bb_\alpha + P^r_i \theta^
i \bb_r $. Substituting this into $\iota_{\be_\alpha} \Omega = 0$ in $\hat 4]$ the coefficients of $\bb_r$ give
$$
P^r_i b_{\alpha j}^i - \hat K_{\alpha s}^r  P ^s_j = 0 .
\tag{A.3}
$$
The coefficients of $\bb_\beta$ in $\iota_{\be_\alpha} \Omega = 0$ give
$$
\iota_{\be_\alpha} d \omega^\beta + K^\alpha_{\beta \gamma} \omega^\gamma +
 \tilde K^\beta_{\alpha r}P^r_i \theta^i = 0 \ .
$$
This equation leads to the formula
$$
d\omega^\alpha = -\textstyle{\frac{1}{2} K^\alpha_{\beta \gamma} \omega^\beta \wedge \omega^ \gamma -\tilde K^\alpha_{\beta r} P^r_i \omega^\beta \wedge \theta^i
 -\frac{1}{2} C^\alpha _{jk} \theta^j\wedge \theta^k}
\tag{A.4}
$$
where $C^\alpha _{jk}=C^\alpha_{[jk]}$ are yet to be determined. The form $\Omega$ can be computed from (A.4) and (A.3):
$$
\Omega =\textstyle{
\left( P^r_i P^s_k \tilde b^i_{sj} +\frac{1}{2}   L^r_{st} P^s_j P^t_k \right) \theta^j \wedge \theta^k\otimes \bb_r +
\frac{1}{2}\left(\tilde L^\alpha_{rs} P^r_iP^s_j -C^\alpha_{ij}\right) \theta^i \wedge \theta^j \otimes \bb_\alpha \ }.
$$

By choosing  $P^r_i$ and $C^i_{jk}$, we can satisfy the last two equations in $\hat 4]$ (the Bianchi identities).
The Bianchi identities can be imposed by either computing $\Omega$ and  implementing them as written in $\hat 4]$ or imposing $d^2\theta^i=0$ and $d^2 \omega^\alpha = 0$ in (A.2) and (A.4).

We write, in terms of our basis, the condition for the algebra
$\lag$ to be reductive.

\proclaim{Lemma A1} Let $G \to G/H$ with $H$  be a homogeneous space with $H$ connected and where the Lie algebra $\lag$ admits forms $\theta$ and $\omega$ satisfying the conditions in Lemma 3.2. Then $G\to G/H$ is reductive if and only if there exists constants $ r^\alpha_i $ such that
$$
r^\gamma_i K^\alpha_{\beta \gamma} -r^\alpha_i b^i_{\beta j}=\tilde K^\alpha_{\beta r} P^r_i \ .
\tag{A.5}
$$
\endproclaim

\demo{Proof} Since $H$ is connected, the form $(\omega^\beta + r^\beta_i \theta^i)\otimes \be_\beta$ is invariant (or equivariant) if and only if its Lie derivative with respect to $\be _\alpha \in \lah $ is zero.  Therefore, $G/H$ is reductive if and only if there exists $r^\beta_i$ such that
$$
\iota_{\be_\alpha} \left( d \omega ^\beta + r^\beta _i d \theta^i \right)\otimes \be_\beta + \left( \omega ^\beta + r^\beta _i \theta^i \right)\otimes K^\gamma_{\alpha \beta }  \be_\gamma =0 \ .
$$
Expanding this equation out using (A.2), (A.3),  and (A.4) we get equation (A.5).
\qed
\enddemo

\head Thanks. \endhead

The authors would like to thank Ian Anderson for numerous helpful suggestions and for the use of his package {\it Vessiot} available at
$www.math.usu.edu/ \ \tilde{ }fg\_mp $. We also thank Charlie Torre, Stephen Yeung, and Scot Adams for helpful discussions.

\InitializeRef

\heading References \endheading

\enddocument